%%!! Last update 99-06-24 09:34:35
%%The final version submitted to Finance and Stochastics
%%Sent off: June 23, 1999

\input amstex         
\input epsf
\def\prodi{\mathop{{\lower 3pt\hbox{\epsfxsize=7pt\epsfbox{pi.ps}}}}}

\documentstyle{amsppt}

\magnification=1200
\vsize=9.4truein
\voffset=-.2truein
\hsize=6.2truein
\baselineskip=35pt
\lineskip=1pt
\lineskiplimit=0pt
\tolerance=500
\TagsOnRight
\headline={\hss\tenrm\folio\hss}
%\pageno=0
\footline={\hfill}
\nologo

\def\RR{\Bbb R}

\def\pindex{\upsilon}
\def\osc{Osc}

\def\path{function}
\def\paths{functions}
                               
%notation for stock price:
\def\spr{P}
%notation for trading strategy:
\def\trs{\phi} 
%notation for portfolio:
\def\prf{V}
\def\Q{\frak Q}
\def\dim{\nu}
\def\g{F}
\def\dg{\g^{\prime}}

\def\dgcfl{\g_l^{\prime}{\circ}f}

%\catcode`@=11
%\def\myblock{\RIfMIfI@\nondmatherr@\block\fi
%       \else\ifvmode\medskip\noindent\fi
%        $$\def\endblock{\par\egroup$$}\fi
%  \vbox\bgroup\advance\hsize-1\indenti\noindent}
%\catcode`@=13
%\catcode`@=11
%\def\block{\RIfMIfI@\nondmatherr@\block\fi
%       \else\ifvmode\vskip\abovedisplayskip\noindent\fi
%        $$\def\endblock{\par\egroup$$}\fi
%  \vbox\bgroup\advance\hsize-2\indenti\noindent}
%\catcode`@=13

\document

\vskip 4pc
\flushpar
MODELLING OF STOCK PRICE CHANGES:
\flushpar
A REAL ANALYSIS APPROACH
\vskip 1pc

\flushpar
{\bf Rimas Norvai\v sa\footnote"*"{Research supported in 
part by an NSERC Canada Collaborative Grant at Carleton University, Ottawa, 
Canada, by a U.S. National Science Foundation Grant
and by an NWO (Dutch Science Foundation) grant.
Manuscript received: ...; final version received: ...}}
\vskip 1pc

\flushpar
{Institute of Mathematics and Informatics,
Akademijos 4, LT-2600 Vilnius, Lithuania
(e-mail: norvaisa\@ktl.mii.lt)} 

\vskip 2pc
\flushpar
{\bf Abstract.} In this paper a real analysis approach to
stock price modelling is considered.
A stock price and its return are defined in a duality to each other
provided there exist suitable limits along a sequence of nested partitions
of a time interval, mimicking sum and product integrals.
It extends the class of stochastic processes susceptible to theoretical
analysis.
Also, it is shown that extended classical calculus is 
applicable to market analysis whenever the local $2$--variation
of sample {\paths} of the return is zero, or is determined by jumps
if the process is discontinuous. 
In particular, an extended Riemann-Stieltjes integral is used in that 
case to prove several properties of trading strategies.

\vskip 2pc
\flushpar
\bigskip
\flushpar
\bigskip
\flushpar

\vskip 2pc

\flushpar
{\bf 1 Introduction and discussion}
\vskip 1pc

\flushpar
In continuous--time financial mathematics the solution to the 
Dol\'eans--Dade stochastic differential equation is often used as a 
model for stock price changes.
The semimartingale driving this equation is called the return.
Since many conclusions on the price behavior depend on the return,
it plays an important role in mathematics of finance.
On the other hand, the returns in econometrics of financial markets
are sometimes modelled by stochastic processes which are not
semimartingales.
To provide a theoretical justification for such cases, one introduces 
a Dol\'eans--Dade type equation with the stochastic integral replaced 
by a different integral.
However, solutions to integral equations based on different
integrals may differ considerably as demonstrated Wong and
Zakai (1965).
One may ask then whether it is possible to build up a model of
stock price changes which is independent of a particular
integration theory?
The present paper addresses this question and provides a new insight 
into the relation between theoretical and applied financial mathematics.

{\sl 1.1. Prices and returns}.
To begin with we discuss two continuous--time
stochastic models for a frictionless stock market.
Let $R=\{R(t)\colon\,0\leq t\leq T\}$ be a semimartingale
such that $R(0)=0$ almost surely and let $Q=\{Q(t)\colon\,0\leq t\leq 
T\}$,  where $Q(t):=\exp\{R(t)\}$ for $0\leq t\leq T$.
The pair $(Q,R)$ will be called the {\it exponential system}
of a stock.
Then $Q$ is the price and $R$ is the return of a stock of
the exponential system $(Q,R)$.
Let $P=\{P(t)\colon\,0\leq t\leq T\}$ be a stochastic process
satisfying the equation
$$
P(t)=1+(SI)\int_0^tP(s-)\,dR(s),\quad 0\leq t\leq T,\tag1.1
$$
where $P(0-):=1$ and $(SI)$ denotes the stochastic integral
defined by the $L^2$-isometry.
Dol\'eans--Dade (1970) proved that the unique solution to (1.1)
is given by 
$$
P(t)=\exp \{R(t)-\frac 12 \langle R^c,R^c\rangle (t)\}\prod_{(0,t]}
(1+\Delta R)\exp\{-\Delta R\},\quad 0< t\leq T,
$$    
and $P(0)=1$,
where $R^c$ is the continuous local martingale part of the
semimartingale $R$ and $\Delta R(s):=R(s)-R(s-)$ for $s\in (0,T]$.
If $P$ satisfies (1.1) and is bounded away from zero then,
by associativity of 
%(the substitution rule for) 
the stochastic integral, we have
$$
R(t)=(SI)\int_0^t\frac {dP(s)}{P(s-)},\quad 0\leq t\leq T.\tag1.2
$$
The pair $(P,R)$ satisfying (1.1) will be called the {\it stochastic 
exponential system}.
Then $P$ is the price and $R$ is the return of a stock of
the stochastic exponential system.
Parts of continuous--time financial mathematics based on the
exponential system and on the stochastic exponential system will be called
respectively the exponential model and the stochastic exponential model.
In general, the exponential system is different from the
stochastic exponential system.
Indeed, if $R$ is a standard Brownian motion $B=\{B(t)\colon\,t\geq 0\}$, 
then the solution to (1.1) is the stochastic process 
$P_B(t):=\exp\{B(t)-t/2\}$, $0\leq t\leq T$, often called
the geometric Brownian motion.
In both systems the prices are observable quantities 
meaning that they represent real data, 
while the returns are non--observable and depend on the models.
In addition to being a semimartingale, $R$ may sometimes be assumed to 
satisfy certain probabilistic conditions about its distribution.
An adequacy to real data of such assumptions on $R$ can be tested 
by using the  price transformations: the log return $R(t)=\log Q(t)$,  
$0\leq t\leq T$, for the exponential model, and the return (1.2) for the 
stochastic exponential model.
The log return is often used in econometric literature
which means that certain hypotheses about the exponential
model are tested.
If one wishes to test the stochastic exponential model then the
return (1.2) has to be used.
However (1.2) is {\it not} defined for a single sample {\path},
so that its statistical tractability is problematic.
On the other hand, as pointed out B\"uhlmann, Delbaen, Embrechts
and Shiryaev (1996),
under probabilistic price analysis, the stochastic exponential model
turns out to be more advantageous than the exponential model.
Therefore, it is appealing to modify the stochastic exponential model
in such a way that to make it more manageable for statistical
analysis.
B\"uhlmann et al. (1996) provide the analysis of the exponential 
model via its reduction to the stochastic exponential 
model using a suitable transformation in (1.1) instead of $R$.

{\sl 1.2. Price changes as an evolution}.
In the present paper we define a price and its return in a duality
without a priori integrability or probabilistic assumptions 
(Definition 2.9 below), 
and show that almost all sample {\paths} of 
many typical stochastic processes including a Brownian motion
satisfy the new definition (Propositions 2.10 and 2.11 below).
An idea behind the definition is based on known results about a
one--to--one correspondence between an evolution and its generator.
A family of real numbers $U=\{U(s,t)\colon\,a\leq s\leq t\leq b\}$ 
is an evolution on $[a,b]$ if $U(t,t)=1$  and
$$
U(t,r)U(r,s)=U(t,s)\quad\text{for all}\quad a\leq s\leq r\leq t\leq b.
$$
Let $P$ be a function representing a stock price over a period
$[0,T]$ such that $P(0)=1$,
and let $U(t,s):=P(t)/P(s)$ for $0\leq s\leq t\leq T$.
Then $U$ so defined is a simple example of
an evolution on $[0,T]$ defined by stock price changes.
An evolution arise in describing the state of {\it nonautonomous}
systems and they are generalizations of the concept of a 
one--parameter semigroup of bounded linear operators on a
Banach space describing the state of autonomous linear systems.
The classical Hille--Yosida theorem describes any strongly
continuous, contractive semigroup in terms of its generator.
In this way the Hille--Yosida theorem provides a one--to--one
correspondence between semigroups and their generators.
An important difficult question is when and in what sense
will a given evolution $U$ have a generator?
The answer depends on a behaviour of the function
$[a,b]\ni t\mapsto U(t,a)$, in particular on its  
$p$-variation.
If $U$ is defined by stock price changes, that is if $U(\cdot,0)=P$, 
then its generator is a return as defined in the present paper.
Therefore the pair $(P,R)$ satisfying Definition 2.9 is called
the ({\it weak}) {\it evolutionary system}.

{\sl 1.3. The $p$-variation}.
For the approach advocated in the present paper,
the notion of $p$-variation of a function plays a role comparable with
a role of a martingale property in the stochastic exponential model.
For a function $f\colon\,[a,b]\mapsto\RR$ and a real number
$0<p<\infty$, the {\it $p$-variation} $v_p(f)=v_p(f;[a,b])$ is the least
upper bound of sums $s_p(f;\kappa):=\sum_{i=1}^n|f(x_i)-f(x_{i-1})|^p$
over all partitions $\kappa=\{x_i\colon\,i=0,\dots,n\}$ of $[a,b]$.
We notice that the $2$-variation is {\it not} the same as the quadratic
variation.
For a standard Brownian motion $B=\{B(t)\colon\,t\geq 0\}$ and any 
$0<T<\infty$, $v_2(B;[0,T])=+\infty$ almost surely, while $v_p(B;[0,T])
<\infty$ for each $p>2$ and the quadratic variation of $B$ is 
defined in the almost sure sense for certain sequences of partitions.
For any function $f$ on $[a,b]$, define the {\it index of $p$-variation}
$\pindex (f)=\pindex (f;[a,b])$ by
$$
\pindex (f;[a,b]):=\cases \inf\{p>0\colon\,v_p(f)<\infty\} &\text{if
the set is nonempty}\\ +\infty&\text{otherwise.}\endcases
$$ 
Therefore for a Brownian motion $B$, $\pindex (B;[0,T])=2$
almost surely.
Also for any $0<T<\infty$, $\pindex (X;[0,T])<2$ almost surely 
if $X$ is a mean zero Gaussian stochastic process with stationary
increments, continuous in quadratic mean and the incremental
variance $\{E[X(t+u)-X(t)]^2\}^{1/2}$ varies regularly as 
$u\downarrow 0$ with index $\gamma > 1/2$, 
or if  $X$ is a homogeneous L\'evy process
with the L\'evy measure $L$ such that $\smallint_{\RR\setminus\{0\}}
(1\wedge |x|^p)\,L(dx)<\infty$ for some $p<2$.

{\sl 1.4. Stochastic and classical calculi}.
In this paper it is proved that evolutionary
systems $(P,R)$ possess a uniformity property with respect to
sequences of partitions defining $P$ and $R$ provided the $p$-variation
index $\pindex (R)<2$.
This fact is important when one deals with fitting a model to
real data, or when one considers a relation between discrete--time
and continuous--time models (cf. Theorem 3.5 below).
With the help of the result of F\"ollmer (1981) one can show
that a weak evolutionary system $(P,R)$ satisfies an integral equation
similar to (1.1) with a different integral.
If sample functions of the return $R$ in the evolutionary system $(P,R)$ 
have the $p$-variation index $\pindex(R)<2$, then (1.1) and (1.2) hold
path by path with the stochastic integral replaced by the Left Young 
integral, an extended Riemann--Stieltjes integral.
Norvai\v sa (1999) proved that the values of the Left Young integral and
the values of the corresponding stochastic integral agree almost surely
under conditions ensuring the existence of both.
In this sense the value $\pindex (R)=2$ of the $p$-variation
index is a borderline between an area where classical calculus applies
and an area where stochastic calculus is needed essentially.
We notice that the semimartingale property of the return $R$
in the stochastic exponential system $(P,R)$ makes a borderline between
classical and stochastic calculi on a different level, e.g.
the value $\pindex (R)=1$ (the $1$-variation is the same as the total
variation).
Several examples of returns such as a hyperbolic L\'evy motion, 
a normal inverse Gaussian L\'evy process, the V.G. process, 
an $\alpha$-stable L\'evy motion with $\alpha\in [1,2)$, 
or a fractional Brownian motion with the Hurst index $H\in (1/2,1)$,
can be treated using classical calculus.
The aim of the present paper is to find a connection between
the two calculi for the mathematics of finance.
However, more interesting is a question whether it is possible to
develop a full fledged model of a financial market based on
the evolutionary system.
Clearly it is not possible to answer to this question
at this writing.
A model construction requires a more advanced development of theories 
of integral equations and optimal control for functions of bounded 
$p$-variation, as well as further development of concepts of market 
efficiency, equilibrium and risk in the new context.
 
{\sl 1.5. Arbitrage}.
We finish with a discussion of arbitrage for the evolutionary system.
In the continuous--time financial mathematics based on the
semimartingale theory the first and second fundamental theorems 
deal with the key principals of the theory.
These theorems relate suitable forms of an arbitrage with the existence 
of a (unique) martingale measure and completeness.
Thus an applicability of these tools is restricted if arbitrage
is possible.
In particular, this concerns the contingent claim valuation theory
based on the no arbitrage principle.
Less formally, the no arbitrage principle is considered as a natural 
property of a model of an {\it ideal} financial market 
because ``there is no such thing as a free lunch'' in equilibria market.
These arguments may give an impression that no approaches other than 
martingale based  stochastic calculus can be useful for mathematical 
finance.
However the situation is not as simple as it may look.
There are examples from a game theory
where such a thing as a free lunch is possible under equilibria
(see p.\ 137 in Kac, Rota and Schwartz, 1992).
On the other hand, non--equilibrium can explain stylized facts discovered 
through the statistical analysis of market data (see Chapter 4
of an overview of Focardi and Jonas, 1997,  based on interviews
with over 100 persons in industry and academia).
A strong critique of a whole current financial mathematics comes from 
actuaries who use different principles to value contingent claims
(see e.g.\ Clarkson, 1996, 1997).
So instead of avoiding arbitrage it seems more fruitful to
have a model which accommodates both, free lunch areas 
as well as areas without a free lunch, and leave the question
of performance evaluation of such a model to econometrics.

Next we illustrate how a real analysis approach may shed  
new light on arbitrage.
One way to define an arbitrage for evolutionary systems
is to follow the pattern from the stochastic exponential model which 
requires first to define a self--financing strategy.
As pointed out Harrison and Pliska (1981, Section 7),
the restriction to predictable trading strategies as well as to
gains defined using the stochastic integral needs a careful study.
Clearly we cannot use these constructions in the present setting.
Instead we define self--financing strategies pathwise following the logic
of the present approach (Definition 3.3 below), 
and prove that the criteria suggested by Harrison and Pliska (1981) 
does apply to the new notions (Theorem 3.5 below).
Then arbitrage can be defined either for a single function representing
a price evolution, or for almost every sample function of a stochastic 
process using notation of Section 3 as follows:
given a price $P=(P_0,\dots,P_{\dim})$ (of $1+\dim$ assets) 
during a time period $[0,T]$, a self--financing $P$--trading strategy 
$\phi=(\phi_0,\dots,\phi_{\dim})$ is an {\it arbitrage opportunity} for 
$P$ at time $T$ if the portfolio value function $V^{\phi,P}$ is $0$ at $0$ 
and positive at $T$.
Salopek (1998) proved that an arbitrage in this sense can be constructed
whenever the return of an evolutionary system is {\it continuous}
function of bounded $p$-variation for some $1\leq p<2$.
To give a short proof of the same statement we modify the ingenious 
construction of an arbitrage due to Shiryaev  (1998, Example VII.2c.4).
To this aim we replace a fractional Brownian motion with the Hurst index 
$H\in (1/2,1)$ in his construction, with a continuous function of bounded
$p$-variation for some $1\leq p<2$, and use the chain rule formula
given by Theorem 2.1 below instead of It\^o's formula.

\proclaim{Proposition 1.1}
Let $f$ be a continuous function of bounded $p$-variation on $[0,T]$
for some $1\leq p<2$ such that $f(0)=0$ and $f(T)\not =0$,
and let $r$, $\sigma$ be real numbers.
Then for the evolutionary system $(P,R)$ with $R(t)=(rt,rt+\sigma 
f(t))$, $0\leq t\leq T$, there exists an arbitrage opportunity
for $P$ at time $T$.
\endproclaim

\demo{Proof}
Let $R_0(t):=rt$ and $R_1(t):=rt+\sigma f(t)$ for $0\leq t\leq T$.
By Proposition 2.6, $P_0(t)=e^{rt}$ and $P_1(t)=e^{rt+\sigma f(t)}$
for $0\leq t\leq T$.
The vector function $P=(P_0,P_1)$ is the price in the sense defined
in Section 3 below.
Let $\phi_0(t):=1-\exp\{2\sigma f(t)\}$ and $\phi_1(t):=
2[\exp\{\sigma f(t)\}-1]$ for $0\leq t\leq T$.
By Proposition 3.2, $\phi=(\phi_0,\phi_1)$ is the $P$--trading strategy.
Next we show that $\phi$ is self--financing $P$--trading strategy as
defined in Definition 3.3.
To this aim we apply the chain rule formula from Theorem 2.1 twice:
first take $h\equiv 1$, $F(u_1,u_2)=u_1u_2^2$,
and second take $h=\phi_1$, $F(u_1,u_2)=u_1u_2$.
For each $0<t\leq T$, we then have
$$
\align
V^{\phi}(t) &=e^{rt}\Big [e^{\sigma f(t)}-1\Big ]^2
=\int_0^t\Big [e^{\sigma f(t)}-1\Big ]^2\,de^{rt}
+\int_0^t2e^{rt}\Big [e^{\sigma f(t)}-1\Big ]\,de^{\sigma f(t)}\\
&=\int_0^t\Big [1-e^{2\sigma f(t)}\Big ]\,de^{rt}
+\int_0^t2\Big [e^{\sigma f(t)}-1\Big ]\,de^{rt+\sigma f(t)}
=V^{\phi}(0)+G^{\phi}(t),
\endalign
$$
where all integrals exist in the Riemann-Stieltjes sense
by the Stieltjes integrability theorem of L.C.\ Young (1936,
p.\ 264).
Since $V^{\phi}(0)=0$ and $V^{\phi}(T)>0$, 
the self--financing $P$--trading strategy $\phi$ is
an arbitrage opportunity for $P$ at time $T$.
\qed\enddemo

The preceding fact shows the irrelevance of a ``long memory''
of a fractional Brownian motion with respect to an arbitrage.
Sample  function behavior of a stochastic process is
responsible for arbitrage opportunities.
A Weierstrass function is a non--probabilistic example of
a function $f$ satisfying hypotheses of Proposition 1.1.
Once we except the evolutionary system as a base for a model
of a financial market then an arbitrage is a property of a local
behavior of a sample function rather than a correlation property
between indefinitely increasing time moments.

\vskip 1pc
\flushpar
{\bf 2 Duality between price and return}
\vskip 1pc

\flushpar
{\sl 2.1. Returns in discrete}--{\sl time models}.
If time $t$ is discrete, say $t=0,1,\dots,T$, there are at least
two different notions of return.
Let $P=\{P(t)\colon\,t=0,1,\dots,T\}$ be a price of
a stock which pays no dividends.
The {\it simple net return} $R_1=\{R_1(t)\colon\,t=0,1,\dots,T\}$
is defined by setting $R_1(0):=0$, and for each $t=1,\dots,T$,
$$
\widehat R_1(t):=R_1(t)-R_1(t-1):=
\cases [P(t)-P(t-1)]/P(t-1),&\text{if $P(t-1)>0$},\\
0,&\text{if $P(t-1)=0$}.\endcases\tag 2.1
$$
Notice that $\widehat R_1(t)$ depends on values $P(t-1)$ and $P(t)$,
so that $\widehat R_1$ is the function of a subinterval $[t-1,t]$.
Notation $\widehat R_1$ (as well as $\widehat R_2$ defined below)
is natural to use in discrete--time models
where time lags have {\it fixed} length.
A work with continuous--time models requires to treat returns
either as interval functions defined on {\it all} subintervals of
$[0,T]$, or as point functions on $[0,T]$.
In this paper we choose to use the form of a point function.
Given $P(0)>0$,
there is a one--to--one correspondence between a positive price $P$ 
and a simple net return $R_1$ having jumps bigger than minus one, 
as described by Pliska (1997, Section 3.2).  
Namely, in addition to (2.1), for each $t=1,\dots,T$, we have
$$
P(t)=P(0)+\sum_{s=1}^tP(s-1)\widehat R_1(s)\quad\text{
and }\quad
P(t)=P(0)\prod_{s=1}^t[1+\widehat R_1(s)].\tag2.2
$$
This correspondence is used in security market models by specifying
simple net returns rather than prices.

Another type of a return is the {\it log return} 
$R_2=\{R_2(t)\colon\,t=0,1,\dots,T\}$ 
defined by setting $R_2(0):=0$ and, for all $t=1,\dots,T$,
$$ 
\widehat R_2(t):=R_2(t)-R_2(t-1):=
\cases \log [P(t)/P(t-1)],&\text{if $P(t-1)>0$},\\
0,&\text{if $P(t-1)=0$}.\endcases
$$
This return is often used in
the econometrics literature on security markets.
It is easy to see that $R_2$ satisfies the additivity property
$$
\widehat R_2(t)+\widehat R_2(t-1)+\dots +\widehat R_2(t-s+1)
=\log [P(t)/P(t-s)]\tag 2.3
$$
for any $s, t\in\{0,1,\dots,T\}$, $s<t$, because the right side 
is the log return for the time period between $t-s$ and $t$.
The additivity property of log returns is one reason of its
popularity among econometricians.
For a discussion of these and other related properties of
returns, see Campbell, Lo and MacKinlay (1997, Section 1.4.1).

Sometimes statistical conclusions based on the log return $R_2$ are
applied to the model (2.2) or even to the continuous--time
stochastic exponential model.
To justify this one usually argues that $\widehat R_1$ and $\widehat R_2$ 
are relatively close to each other when both are small.
However, there are cases when the difference between
$\widehat R_1$ and $\widehat R_2$ cannot be neglected (see e.g.
Elton, Gruber and Kleindorfer, 1975).
To test the stochastic exponential model one needs to use the return $R$
defined via the It\^o integral (1.2).
Then one has to be able to evaluate the It\^o integral using finitely
many values of a single sample {\path}. 
Given a sequence of partitions into shrinking subintervals of $[0,T]$,
by the dominated convergence in probability theorem, the value of $R$
can be approximated in probability by corresponding 
Riemann--Stieltjes sums.
However, it is not possible to conclude the convergence with probability 1
without further restrictions.
A relationship between the price and its return is
suggested below for the continuous--time framework which makes this
approximation possible path by path.
This relationship is motivated by a duality relation between 
additive and multiplicative interval functions, which is also known
as the evolution representation problem.
Recall that the additivity property is satisfied by the log 
return $R_2$ (see (2.3)) while the simple net return $R_1$ lacks this 
property.
On the other hand, the multiplicativity property
$$
\widehat P(t)\cdot\widehat P(t-1)\cdot\dots\cdot\widehat P(t-s+1)=P(t)/P(t-s)
$$
for any $s,t\in\{0,1,\dots,T\}$, $s<t$, is satisfied by the price ratios
$\widehat P(t):=P(t)/P(t-1)$, $t=1,\dots,T$.

{\sl 2.2. The chain rule formula}.
For a finite interval $J$, open or closed at either end,
let $Q(J)$ be the set of all partitions
$\kappa=\{x_i\colon\,i=0,\dots,n\}$ of $J$.
As before for $f\colon\,J\mapsto\RR$ and $0<p<\infty$, let
$v_p(f):=v_p(f;J):=\sup\{s_p(f;\kappa)\colon\,\kappa\in Q(J)\}$
be the $p$-variation of $f$,
where $s_p(f;\kappa):=\sum_{j=1}^n|f(x_i)-f(x_{i-1})|^p$ for
$\kappa=\{x_i\colon\,i=0,\dots,n\}$.
Denote by ${\Cal W}_p={\Cal W}_p(J)$ the set of all functions
$f$ such that $v_p(f)<\infty$.
If $f\in {\Cal W}_p$ for some $p<\infty$ then $f$ is {\it regulated},
that is there exist the limits $f(x-):=\lim_{y\uparrow x}f(y)$ and 
$f(x+):=\lim_{y\downarrow x}f(y)$ when these are defined.
The class of all regulated functions on $J$ will be denoted by
${\Cal R}(J)$.
Given a regulated function $f$ on $[a,b]$, define a left-continuous 
function $f_{-}^{(a)}$ and a right-continuous function $f_{+}^{(b)}$ by
$$
\left\{\aligned
f_{-}^{(a)}(x)&:=f_{-}(x):=f(x-)\quad\text{for}\quad
a<x\leq b\quad\text{and}\quad f_{-}^{(a)}(a):=f(a)\\
f_{+}^{(b)}(x)&:=f_{+}(x):=f(x+)\quad\text{for}\quad
a\leq x<b\quad\text{and}\quad f_{+}^{(b)}(b):=f(b).
\endaligned\right.
$$
Given a regulated function $f$ on $J$, define $\Delta^{-}f$ on $J$ by
$\Delta^{-}f(x):=f(x)-f(x-)$ if $f(x-)$ is defined and $\Delta^{-}f(x)
:=0$ otherwise.
Similarly define $\Delta^{+}f$ on $J$ by $\Delta^{+}f(x):=f(x+)-f(x)$
if $f(x+)$ is defined and $\Delta^{+}f(x):=0$ otherwise.
Since each regulated function $f$ has at most countably many jumps, 
one can define
$\frak S_p (f):=\frak S_p (f;J)
:=\{\sum_{J}(|\Delta^{-}f|^p+|\Delta^{+}f|^p)\}^{1/p}$.
The {\it local $p$-variation} $v_p(f)^{\ast}:=v_p(f;J)^{\ast}$
is defined by
$$
v_p(f;J)^{\ast}:=\inf_{\lambda\in Q(J)}\sup\{s_p(f;\kappa)\colon\,
\lambda\subset\kappa\in Q(J)\}.
$$
Then we have the relation
$\frak S_p (f)^p\leq v_p^{\ast}(f)\leq v_p(f)$.
Let $\Cal W_p^{\ast}=\Cal W_p^{\ast}(J)
:=\{f\in\Cal W_p\colon\,\frak S_p (f)^p=v_p^{\ast}(f)\}$
for $1< p<\infty$.
%Let
%$$
%\Cal D(I):=\{f\in\Cal R(I)\colon\,\Delta^{-}f\Delta^{+}f=0\}.\tag A.2
%$$
For regulated functions $h$ and $f$ on $[a,b]$,
define the {\it Left Young integral}, or the $LY$ integral,
by
$$
(LY)\int_a^bh\,df:=(RS)\int_a^bh_{-}^{(a)}\,df_{+}^{(b)}
+\big [h\Delta^{+}f\big ](a)+\sum_{(a,b)}\Delta^{-}h\Delta^{+}f
$$
provided the Riemann--Stieltjes integral exists in the refinement
sense and the sum  converges absolutely.
Additivity on adjacent intervals as well as some other properties
of the $LY$ integral are proved in Norvai\v sa (1999).
From the Stieltjes integrability theorem of L.C. Young (1936) it 
follows that $(LY)\smallint_a^bh\,df$ is defined if 
$h\in {\Cal W}_p$, $f\in {\Cal W}_q$ and $1/p+1/q>1$.
The following theorem of Norvai\v sa (1999) extends this
result to the case when $1/p+1/q=1$ under additional assumptions
on $h$ and $f$.
Let $\dim$ be a positive integer, and let $F$ be a real-valued
function defined on an open set $U\subset\RR^{\dim}$ containing
a $\dim$-dimensional cube $[c,d]^{\dim}:=[c,d]\times
\cdots\times [c,d]$.
We write $\g\in\Lambda_{1,\alpha}([c,d]^{\dim})$ for $\alpha
\in (0,1]$ if $F$ is differentiable on $U$ with partial 
derivatives $F'_l$, $l=1,\dots,{\dim}$,  and 
there is a finite constant $K_{\alpha}$ such 
that the inequality 
$$
\max_{1\leq l\leq \dim}|\dg_l(u)-\dg_l(v)|\leq
K_{\alpha}\sum_{k=1}^\dim |u_k-v_k|^{\alpha}
$$
holds for all $u=(u_1,\dots,u_\dim), v=(v_1,\dots,v_\dim)\in 
[c,d]^{\dim}$.

\proclaim{Theorem 2.1}
For $\alpha\in (0,1]$, let $f=(f_1,\dots,f_\dim)\colon\,[a,b]\mapsto 
(c,d)^{\dim}$ be a vector function with coordinate functions 
$f_l\in\Cal W_{1+\alpha}^{\ast}([a,b])$ for $l=1,\dots,\dim$,
let $\g\in\Lambda_{1,\alpha}([c,d]^{\dim})$ and
let $h$ be a regulated function on $[a,b]$.
Then the equality
$$
(LY)\int_a^bh\,d(F{\circ}f)=\sum_{l=1}^{\dim} (LY)\int_a^bh(\dgcfl)\,df_l
$$
$$
+\sum_{(a,b]}h_{-}\big [\Delta^{-} (F{\circ}f)-\sum_{l=1}^{\dim}
(\dgcfl)_{-}\Delta^{-}f_l\big ]
+\sum_{[a,b)}h\big [\Delta^{+} (F{\circ}f)-\sum_{l=1}^{\dim}
(\dgcfl)\Delta^{+}f_l\big ]
$$
holds meaning that all ${\dim}+1$ integrals exist provided any $d$ 
integrals exist, and the two sums converge absolutely.
\endproclaim

We refer to the preceding statement as the chain rule formula.
Its proof is given by Norvai\v sa (1999).

{\sl 2.3. Duality relation}.
Turning to a continuous--time model, consider an interval $[0,T]$,
$0<T<\infty$.
Roughly speaking, to extend (2.1) and (2.2) to functions 
defined on $[0,T]$, we pass to a limit along a nested sequence $\lambda$ of 
partitions of $[0,T]$.
The first limit $\Cal L_{\lambda}(f)$ if exists is an extension of the 
(sum) integral,
and for $f\in\Cal W_2^{\ast}$, its values coincide with values of
the $LY$ integral.
The second limit $\Cal E_{\lambda}(g)$ if exists is an extension
of the product integral, and for sample {\paths} $g$ of a semimartingale,
its values coincide with values of the solution to
the Dol\'eans-Dade equation (1.1).
First we prove the existence of $\Cal L_{\lambda}(f)$ and 
$\Cal E_{\lambda}(g)$ for functions $f$ and $g$ from the class 
$\Cal W_2^{\ast}$.
Then the duality relations (2.13) are derived for such functions.
Finally, a duality relation is proved for functions having defined
the quadratic variation.

\definition{Definition 2.2}
Let $\Q([0,T])$ be the set of all nested sequences 
$\lambda=\{\lambda (m)\colon\,m\geq 1\}$ of partitions
$\lambda (m)=\{0=t_0^m<\cdots <t_{n(m)}^m=T\}$ of $[0,T]$  such that 
$\cup_m\lambda (m)$ is dense in $[0,T]$.
Let $I_T$ be either $[0,T]$ or $[0,T)$, and
let $f$ be a real--valued function on $I_T$.
Given $\lambda\in\Q ([0,T])$, we say that {\it $\Cal L_{\lambda}
=\Cal L_{\lambda}(f)$ is defined} on $I_T$ 
if the limit
$$
\Cal L_{\lambda}(f)(t):=\lim_{m\to\infty}\sum_{i=1}^{n(m)}
[f(t_i^m\wedge t)-f(t_{i-1}^m\wedge t)]/f(t_{i-1}^m\wedge t)\tag 2.4
$$
exists for each $t\in I_T$.
Given a nonempty subset $\Q\subset\Q([0,T])$, 
if $\Cal L_{\lambda}(f)$ is defined for 
and does not depend on each $\lambda\in\Q$, then set
$\Cal L_{\Q}=\Cal L_{\Q}(f)$ to be equal to any $\Cal L_{\lambda}(f)$,
$\lambda\in\Q$.
\enddefinition

For a regulated function $f$, a typical example of $\Q\subset\Q([0,T])$ 
in the preceding definition is the set $\Q(f)$ defined by
$$
\Q(f):=\cases \Q([0,T]),&\text{if $f\in\Cal D(I_T)$,}\\
\big\{\lambda\in\Q([0,T])\colon\,\cup_m\lambda (m)\supset
\Delta_f(I_T)\big\},&\text{if $f\in\Cal R(I_T)\setminus
\Cal D(I_T)$,}\endcases
\tag2.5
$$
where $f\in\Cal D(I_T)$ if, at each point of $(0,T)$, $f$ is either
right-continuous or left-continuous and
$\Delta_f(I_T):=\{x\in (0,T)\colon\,\Delta^{-}f(x)\neq 0$ or 
$\Delta^{+}f(x)\neq 0\}$.
Next, under stated conditions we show that $\Cal L_{\Q(f)}$ is defined
and has values of the indefinite $LY$ integral. 

\proclaim{Proposition 2.3}
Let $f\in\Cal W_2^{\ast}(I_T)$ and let $\inf\{f(t)\colon\,t\in I_T\}
\geq\delta$ for some $\delta >0$.
Then $\Cal L_{\Q(f)}(f)$ is defined on $I_T$.
Moreover, for each $t\in I_T$, $f^{-1}$ is $LY$ integrable
with respect to $f$ on $[0,t]$ and the relation
$$
\Cal L_{\Q(f)}(f)(t)=(LY)\int_0^t\frac {df}{f}
=\log\frac {f(t)}{f(0)}
-\sum_{(0,t]}\left [\log\frac {f}{f_{-}}-\frac {\Delta^{-}f}{f_{-}}
\right ]-\sum_{[0,t)}\left [\log\frac {f_{+}}{f}-\frac {\Delta^{+}f}
{f}\right ]\tag 2.6
$$
holds, where the two sums
converge absolutely.
\endproclaim

For the proof we need an auxiliary statement,
where $Q(S):=\{\kappa\in Q(J)\colon\,\kappa\subset S\}$ for
any subset $S\subset J$.

\proclaim{Lemma 2.4}
Let $f\in\Cal W_p^{\ast}([a,b])$ for some $1<p<\infty$, and
let $S$ be a dense subset of $[a,b]$ containing 
all discontinuity points of $f$.
For each $\epsilon >0$, there exists
$\lambda\in Q(S)$ such that $\sum_{j=1}^kv_p(f;(z_{j-1},z_j))<\epsilon$
whenever $\lambda\subset\{z_j\colon\,j=0,\dots,k\}\in Q([a,b])$.
\endproclaim

\demo{Proof}
Let $S\subset [a,b]$ be as in the statement.
Then $v_p^{\ast}(f;[a,b])$ is equal to the greatest lower bound
of sums $\sum_{j=1}^kv_p(f;[z_{j-1},z_j])$ over
$\{z_j\colon\,j=0,\dots,k\}\in Q(S)$.
Let $\epsilon >0$.
Since $\frak S_p(f)<\infty$, there exists a finite set $\mu\subset
[a,b]$ such that $\sum_{\nu}(|\Delta^{-}f|^p+|\Delta^{+}f|^p)>
\frak S_p(f)^p-\epsilon/2$ for each $\nu\supset\mu$.
Then one can choose $\lambda\in Q(S)$ such that $\lambda
\supset\mu$ and $\sum_{j=1}^kv_p(f;[z_{j-1},z_j])<v_p^{\ast}(f)
+\epsilon /2$ whenever $\lambda\subset\{z_j\colon\,j=0,\dots,k\}\in
Q([a,b])$.
For each small enough $\delta >0$ and each $j=1,\dots,k$, we have  
$v_p(f;[z_{j-1},z_j])\geq v_p(f;[z_{j-1},z_{j-1}+\delta])+v_p(f;[z_{j-1}
+\delta,z_j-\delta])+v_p(f;[z_j-\delta,z_j])$.
Letting $\delta\downarrow 0$ and using Lemma 2.19 of 
Dudley and Norvai\v sa (1999), we get
$$
\align
\sum_{j=1}^kv_p(f;(z_{j-1},z_j))&\leq\sum_{j=1}^kv_p(f;[z_{j-1},z_j])
-\sum_{j=1}^k\big [|\Delta^{-}f(z_j)|^p+|\Delta^{+}f(z_{j-1})|^p
\big ]\\
&<V_p^{\ast}(f)^p+\epsilon/2-\frak S_p(f)^p+\epsilon/2=\epsilon.
\endalign
$$
The proof of Lemma 2.4 is complete.
\qed\enddemo

\demo{Proof of Proposition 2.3}
The existence of the integral $(LY)\smallint_0^tf^{-1}\,df$ and
the second equality in (2.6) for each $t\in I_T$ follow from Theorem 
2.1.
To see if it's true take $F(u):=\log u$ for $u\in 
[\delta,\|f\|_{\infty}]$, $\dim =1$, $\alpha =1$, $h\equiv 1$,
and notice that $\Delta^{-}(F{\circ}f)=\log (f/f_{-})$,
$\Delta^{+}(F{\circ}f)=\log (f_{+}/f)$.
To prove that $\Cal L_{\Q(f)}(f)$ is defined on $I_T$ and that the first
equality in (2.6) holds, for each $u\in (0,T]\cap I_T$ and $v\in [0,T)$, let
$$
\phi^{-}(u):=\log\frac {f(u)}{f(u-)}-\frac {\Delta^{-}f(u)}{f(u-)}
\quad\text{and}\quad\phi^{+}(v):=\log\frac {f(v+)}{f(v)}-\frac 
{\Delta^{+}f(v)}{f(v)}.
$$
To begin with the second case in (2.5) consider $\{\lambda (m)\colon\,
m\geq 1\}\in\Q ([0,T])$ such that $\cup_m\lambda (m)$ contains
all discontinuity points of $f$.
For each $m\geq 1$, let
$$
\psi_m(i):=\log\frac {f(t_i^m)}{f(t_{i-1}^m)}-\left [\frac {f(t_i^m)}
{f(t_{i-1}^m)}-1\right ]\qquad\text{for $i=1,\dots,n(m)$.}
$$ 
Let $S:=\cup_m\lambda (m)$ and $\epsilon >0$.
By Lemma 2.4, and because $f$ is regulated
and the two sums in (2.6) converge absolutely,
one can choose $\kappa:=\{z_j\colon\,j=0,\dots,k\}\in Q(S)$
such that
$$
\sum_{j=1}^kv_2(f;(z_{j-1},z_j))<\epsilon,
\quad
\max_{1\leq j\leq k}\osc (f;(z_{j-1},z_j))<\frac {\delta}{2}
\quad\text{and}\quad
\sum_{\nu}\big (|\phi^{-}|+|\phi^{+}|\big )<\epsilon
$$
for any $\nu\subset I_T\setminus\kappa$.
For each partition $\lambda (m)=\{t_i^m\colon\,i=0,\dots,n(m)\}$ 
containing $\kappa$ and for each 
$j\in\{0,\dots,k\}$, let $i(j)\in\{0,\dots,n(m)\}$ be an index such 
that $z_j=t_{i(j)}^m$.
Since
$$
\lim_{m\to\infty}\psi_m(i(j))=\phi^{-}(z_j)\quad\text{and}\quad
\lim_{m\to\infty}\psi_m(i(j-1)+1)=\phi^{+}(z_{j-1})\tag2.7
$$
for $j=1,\dots,k$, one can choose an integer $M\geq 1$ such that
$\lambda (M)\supset\kappa$, there are at least two elements of
$\lambda (M)$ in each interval $(t_{i(j-1)}^m, t_{i(j)}^m)$, 
$j=1,\dots,k$, and
$$
\Big |\sum_{j=1}^k\Big [\psi_m(i(j))-\phi^{-}(z_j)\Big ]+
\sum_{j=0}^{k-1}\Big [\psi_m(i(j)+1)-\phi^{+}(z_j)\Big ]\Big |
<\epsilon
$$
for $m\geq M$.
Let $t\in I_T$.
First suppose $t\in S$, so that $t=t_{l(m)}^m$ for some $l(m)\in
\{1,\dots,n(m)\}$ and for all $m$ larger than some $N(t)$.
For each $m\geq M\vee N(t)$, let $l:=\max\{j\leq k\colon\,z_j\leq t\}$
and $J:=\{0,\dots,l(m)\}\setminus\{i(j),i(j-1)+1
\colon\,j=1,\dots,l\}$.
By the Taylor series expansion with remainder, we have
$|\log (1+u)-u|\leq 2u^2$ for each $|u|\leq 1/2$.
Then, for all $m\geq M\vee N(t)$, we get
$$
\Big |\log\frac {f(t)}{f(0)}-\sum_{(0,t]}\phi^{-}-\sum_{[0,t)}
\phi^{+}-\sum_{i=1}^{n(m)}\big [\frac {f(t_i^m\wedge t)}
{f(t_{i-1}^m\wedge t)}-1\big ]\Big |
$$
$$
<2\epsilon +
2\sum_{i\in J}\left |\frac {f(t_i^m)}{f(t_{i-1}^m)}-1\right |^2
<2\epsilon +\frac {2}{\delta^2}\sum_{j=1}^kv_2(f;(z_{j-1},z_j))
<2\epsilon (1+\delta^{-2}).\tag2.8
$$
If $t\in I_T\setminus S$ then we have in addition the term
$$
\Big |\log\frac {f(t)}{f(t_{l(m)}^m)}-\Big [\frac {f(t)}{f(t_{l(m)}^m)}
-1\Big ]-\sum_{(t_{l(m)}^m,t]}\phi^{-}-\sum_{[t_{l(m)}^m,t)}\phi^{+}
\Big |,
$$
where $l(m):=\max\{i\leq n(m)\colon\,t_i^m<t\}$.
This term tends to zero as $m\to\infty$ because $f$ is continuous at $t$
in this case.
Since $\epsilon$ in (2.8) is arbitrary, $\Cal L_{\Q(f)}(f)$ is defined 
on $I_T$ and the first equality in (2.6) holds for the second case in (2.5).
The proof when $f\in\Cal D(I_T)$ is the same except that we use Lemma 
2.4 with $S=I_T$, and choose $\{t_i(j)^m\colon\,j=0,\dots,k\}$ so that
$z_j\in (t_{i(j)-1}^m,t_{i(j)}^m]$ if $f$ is right--continuous at $z_j$
and $z_j\in [t_{i(j)}^m,t_{i(j)+1}^m)$ if $f$ is left--continuous
at $z_j$.
The proof of Proposition 2.3 is complete.\qed
\enddemo

Using notation as in Definition 2.2, we have:

\definition{Definition 2.5}
Let $g$ be a real--valued function on $I_T$.
Given $\lambda\in\Q ([0,T])$,
we say that {\it $\Cal E_{\lambda}=\Cal E_{\lambda}(g)$ is defined} 
on $I_T$ if the limit
$$
\Cal E_{\lambda}(g)(t):=\lim_{m\to\infty}\prod_{i=1}^{n(m)}
[1+g(t_i^m\wedge t)-g(t_{i-1}^m\wedge t)]\tag2.9
$$
exists for each $t\in I_T$.
Given a nonempty subset $\Q\subset\Q([0,T])$, 
if $\Cal E_{\lambda}(g)$ is defined for 
and does not depend on each $\lambda\in\Q$, then we define
$\Cal E_{\Q}=\Cal E_{\Q}(g)$ to be equal to any $\Cal E_{\lambda}(g)$,
$\lambda\in\Q$.
\enddefinition

Next, under stated conditions we show that (2.9) is defined  
and has values of the {\it product integral with respect to $g$ over}
$[0,t]$,  $\prodi_0^t(1+dg)$,
defined as the limit of the product from $i=1$ to $n$
of $1+g(t_i)-g(t_{i-1})$, if it exists, under refinements of
partitions $\{t_i\colon\,i=0,\dots,n\}$ of $[0,t]$.
The set $\Q(g)$ in the following statement is defined by (2.5).

\proclaim{Proposition 2.6}
Let $g\in\Cal W_2^{\ast}(I_T)$ and let $(\Delta^{-}g)\wedge(\Delta^{+}g)>-1$ 
on $I_T$.
Then $\Cal E_{\Q(g)}(g)$ is defined on $I_T$.
Moreover, for each $t\in I_T$, the product integral $\prodi_0^t (1+dg)$ 
exists, is positive and the relation
$$
\Cal E_{\Q(g)}(g)(t)=\prodi_0^t (1+dg)
=e^{g(t)-g(0)}\prod_{[0,t]}\big [(1+\Delta^{-}g)(1+\Delta^{+}g)\big ]
e^{-\Delta^{-}g-\Delta^{+}g}\tag2.10
$$ 
holds, where the product converges absolutely.
\endproclaim

\demo{Proof}
The product integral  $\prodi_0^t(1+dg)$ exists, and
the second equality in (2.10) holds for each $t\in I_T$ by Theorem 
4.4 of Dudley and Norvai\v sa (1999).
Since all jumps of $g$ are bigger than $-1$, the positivity of the
product integral follows from its definition.
To prove that $\Cal E_{\Q(g)}(g)$ is defined on $I_T$ and 
the first equality in (2.10) holds, let $t\in I_T$ and 
$\{\lambda (m)\colon\,m\geq 1\}\in\Q(g)$, where
$\lambda (m)=\{t_i^m\colon\,i=0,\dots,n(m)\}$.
One can assume that $t=t_{l(m)}^m$ for some $1\leq l(m)\leq n(m)$.
Otherwise we include $t$ into $\lambda (m)$ and change indices.
For a finite set $\mu\subset I_T$, let
$$
A(g;\mu):=\prod_{z\in\mu}\big [(1+
\Delta^{-}g(z))(1+\Delta^{+}g(z)\big ]e^{-\Delta^{-}g(z)
-\Delta^{+}g(z)}.
$$
To begin with the second case in (2.5) consider $\{\lambda (m)\colon\,
m\geq 1\}\in\Q([0,T])$ such that $\cup_m\lambda (m)$ contains all
discontinuity points of $g$.
Let $S:=\cup_m\bar\lambda (m)$, where $\bar\lambda (m)=\{0=t_0^m<\dots
<t_{l(m)}^m=t\}$, and $\epsilon\in (0,2A)$, where
$A(g)$ denotes the product $\prod_{[0,t]}$ in (2.10).
Because $g$ is regulated and by Lemma 2.4,
one can choose $\kappa =\{z_j\colon\,j=0,\dots,k\}\in
Q(S)$ such that $\osc (g(z_{j-1},z_j))<1/2$ for $j=1,\dots,k$,
$$
\sum_{j=1}^kv_2(g;(z_{j-1},z_j))<\frac {\epsilon}{8eA}\quad\text{and}
\quad\big |A(g;\mu)-A(g)\big |<\frac {\epsilon}{4}
$$
for all $\mu\supset\kappa$.
For each $\bar\lambda(m)\supset\kappa$ and for each $j\in\{0,\dots,k\}$,
let $i(j)\in\{0,\dots,l(m)\}$ be such that $z_j=t_{i(j)}^m$.
Let $\Delta_i^mg:=g(t_i^m)-g(t_{i-1}^m)$ for $i=1,\dots,l(m)$, $m\geq 
1$, and let
$$
U(g;\kappa)
:=\prod_{j=1}^k(1+\Delta_{i(j)}^mg)(1+\Delta_{i(j-1)+1}^mg)
\exp\big\{-\Delta_{i(j)}^mg-\Delta_{i(j-1)+1}^mg\big\}.
$$
Then letting 
$J:=\{0,\dots,l(m)\}\setminus\{i(j),i(j-1)+1\colon\,j=1,\dots,k\}$,
we have
$$
\prod_{i=1}^{l(m)}(1+\Delta_i^mg)
=e^{g(t)-g(0)}U(g;\kappa)\exp\big\{\sum_{i\in J}\big [
\log (1+\Delta_i^mg)-\Delta_i^mg\big ]\big\},\tag2.11
$$
for all  $\bar\lambda (m)\supset\kappa$.
Since $|\Delta_i^mg|\leq 1/2$ for $i\in J$, by the Taylor series 
expansion with remainder, we get
$$
\Big |\sum_{i\in J}\big [\log (1+\Delta_i^mg)-\Delta_i^mg\big ]\Big |
=\sum_{i\in J}\theta (\Delta_i^mg)(\Delta_i^mg)^2\leq 2\sum_{j=1}^k
v_2(g;(z_{j-1},z_j))<\frac {\epsilon}{4eA},
$$
where $\theta (u)\in [2/9,2]$ for $|u|\leq 1/2$.
Let $M\geq 1$ be an integer such that $\bar\lambda (M)\supset\kappa$,
there are at least two elements of $\bar\lambda (M)$
in each interval $(t_{i(j-1)}^m, t_{i(j)}^m)$, $j=1,\dots,k$, 
and $|1-U(g;\kappa)/A(g;\kappa)|<\epsilon /(8eA)$ for all $m\geq M$.
Using the inequality $|ue^v-1|\leq e|v|+2e|u-1|$ for $|v|\leq 1/2$
and $|1-u|\leq 1/4$ one can show
that (2.11) differs from the right side of (2.10)
by $\epsilon\exp\{g(t)-g(0)\}$ for all $m\geq M$.
Since $\epsilon$ is arbitrary, $\Cal E_{\Q(g)}(g)$ is defined
on $I_T$ and the first relation in (2.10) holds when $\Q(g)$ is defined
by the second case in (2.5).
The proof for the first case of (2.5) is similar and therefore is omitted.
The proof of Proposition 2.6 is complete.
\qed\enddemo

To show a duality between $\Cal E_{\Q}(g)$ and $\Cal L_{\Q}(f)$ for
$g, f\in\Cal W_2^{\ast}$ first we prove it between the {\it indefinite product
integral} $\Cal P(g)$ and the {\it indefinite $LY$ integral} 
$\Cal S(f)$ defined by
$$
\Cal P(g)(t):=\prodi_0^t(1+dg)\quad\text{and}\quad
\Cal S(f)(t):=(LY)\int_0^t\frac {df}{f}\tag2.12
$$
for $t\in I_T$ whenever the integrals exist.
The following theorem was proved by Dudley and Norvai\v sa (1999,
Theorem 6.8 and 6.10) for functions $f, g$ with values in a Banach algebra 
under the stronger assumption: $f, g\in\Cal W_p$ for some $p\in (0,2)$.
To extend this result to real--valued functions from the class
$\Cal W_2^{\ast}$ we use the chain rule formula.

\proclaim{Theorem 2.7}
{\rm I.} Let $g\in\Cal W_2^{\ast}(I_T)$ and $(\Delta^{-}g)\!\wedge
\!(\Delta^{+}g)>-1$ on $I_T$.
Then the indefinite product integral $\Cal P(g)$ is defined and
$\Cal P(g)\in\Cal W_2^{\ast}(I_T)$.
Moreover, the indefinite $LY$ integral in {\rm (2.12)} is defined for
$f=\Cal P(g)$ and  $\Cal S(\Cal P(g))(t)=g(t)-g(0)$ for $t\in I_T$.

{\rm II.} Let $f\in\Cal W_2^{\ast}(I_T)$ and $\inf\{f(t)\colon\,t\in I_T\}
\geq\delta$ for some $\delta >0$.
Then the indefinite $LY$ integral $\Cal S(f)$ is defined and
$\Cal S(f)\in\Cal W_2^{\ast}(I_T)$.
Moreover, the product integral in {\rm (2.12)} is defined for
$g=\Cal S(f)$ and $\Cal P(\Cal S(f))(t)=f(t)/f(0)$ for $t\in I_T$.
\endproclaim

\demo{Proof}
I. The indefinite product integral $\Cal P(g)$ is defined on $I_T$
by Theorem 4.4 of Dudley and Norvai\v sa (1999).
It is easy to prove that
$\Cal P(g)/\Cal P(g)_{-}=1+\Delta^{-}g>0$, $\Cal P(g)_{+}/
\Cal P(g)=1+\Delta^{+}g>0$ on $I_T$, and
$\Cal P(g)\in\Cal W_2^{\ast}(I_T)$.
Let $t\in I_T$.
By Theorem 2.1, the 
$LY$ integral $\Cal S(\Cal P(g))(t)$ exists and
$$
\Cal S(\Cal P(g))(t)=g(t)-g(0)+\log\Big [\prod_{(0,t]}(1+\Delta^{-}g)
e^{-\Delta^{-}g}\Big ]+\log\Big [\prod_{[0,t)}(1+\Delta^{+}g)
e^{-\Delta^{+}g}\Big ]
$$
$$
-\sum_{(0,t]}\left [\log (1+\Delta^{-}g)-\Delta^{-}g\right ]
-\sum_{[0,t)}\left [\log (1+\Delta^{+}g)-\Delta^{+}g\right ]
=g(t)-g(0).
$$
The last equality follows by taking the limit 
of $\log [\prod_{\mu}\Phi]=\sum_{\mu}[\log\Phi]$
along a nested sequence of finite sets $\mu$ of jump points
of $\Phi=(1+\Delta g)\exp (-\Delta g)$.

II. By Theorem 2.1, $\Cal S(f)$ is defined on $I_T$ and its value is given
by the right side of (2.6).
Since $\log f\in\Cal W_2^{\ast}$ and the two sums in (2.6)
converge absolutely, it follows that $\Cal S(f)\in\Cal W_2^{\ast}$.
Let $t\in I_T$.
By Theorem 4.4 of Dudley and Norvai\v sa (1999) and by Proposition 6
of Norvai\v sa (1999),
the product integral $\Cal P(\Cal S(f))(t)$ exists and has the 
representation
$$
\Cal P(\Cal S(f))(t)=\frac {f(t)}{f(0)}
\exp\Big\{-\sum_{(0,t]}\Big [\log\Big (\frac {f}{f_{-}}\Big )
-\frac {\Delta^{-}f}{f_{-}}\Big ]-\sum_{[0,t)}\Big [\log\Big (
\frac {f_{+}}{f}\Big )-\frac {\Delta^{+}f}{f}\Big ]\Big\}\times
$$
$$
\times\prod_{(0,t]}\Big [(1+\frac {\Delta^{-}f}{f_{-}})\exp (-\frac 
{\Delta^{-}f}{f_{-}} )\Big ]\prod_{[0,t)}\Big [(1+\frac {\Delta^{+}f}{f})
\exp (-\frac {\Delta^{+}f}{f} )\Big ]=\frac {f(t)}{f(0)},
$$
where the last equality follows using the limiting argument as in
the part I.
The proof of Theorem 2.7 is complete.
\qed\enddemo

By Propositions 2.3, 2.6 and Theorem 2.7, for $\Q=\Q([0,T])$ and
for each $t\in I_T$, 
it follows that
$$
\Cal E_{\Q}(\Cal L_{\Q}(f))(t)=f(t)/f(0)\quad\text{and}\quad
\Cal L_{\Q}(\Cal E_{\Q}(g))(t)=g(t)-g(0)\tag2.13
$$
whenever $f, g\in\Cal W_2^{\ast}(I_T)$ are either right-- or 
left--continuous at each point, $f$ is bounded
away from zero and all jumps of $g$ are bigger than $-1$.
Next we partly extend the duality between $\Cal E_{\lambda}$
and $\Cal L_{\lambda}$ for $\lambda\in\Q([0,T])$ and for certain 
functions outside of the class $\Cal W_2^{\ast}$.

\proclaim{Proposition 2.8}
Let $g\in\Cal W_p([0,T])$, $1\leq p<3$, be continuous and let
$\lambda =\{\lambda (m)\colon\,m\geq 1\}$ be a sequence of partitions
$\lambda (m)=\{0=t_0^m<\dots <t_{n(m)}^m=T\}$ such that
the mesh $|\lambda (m)|\to 0$ with $m\to\infty$.

{\rm (I)} For each $t\in [0,T]$, the limit
$$
b_{\lambda}(g)(t):=\lim_{m\to\infty}\sum_{i=1}^{n(m)}
[g(t_i^m\wedge t)-g(t_{i-1}^m\wedge t)]^2\tag2.14
$$
exists if and only if {\rm (2.9)} so does, and then
$$
\Cal E_{\lambda}(g)(t)=\exp\{g(t)-g(0)-2^{-1}b_{\lambda}(g)(t)\}.
\tag2.15
$$

{\rm (II)} Suppose $b_{\lambda}(g)$ from statement {\rm (I)}
is defined and continuous on $[0,T]$.
Then, for each $t\in [0,T]$,  the limit {\rm (2.4)} exists for
$f=\Cal E_{\lambda}(g)$, and
satisfies the relation 
$$
\Cal L_{\lambda}(\Cal E_{\lambda}(g))(t)=g(t)-g(0).\tag2.16
$$
\endproclaim

\demo{Proof}
To prove statement (I), let $t\in (0,T]$.
For each $m\geq 1$, let $l(m)\in\{1,\dots,n(m)\}$ be an integer
such that $t\in (t_{l(m)-1}^m,t_{l(m)}^m]$
and let $u_i^m:=g(t_i^m\wedge t)-g(t_{i-1}^m)$ for $i=1,\dots,l(m)$.
Since $g$ is continuous, there exists
an integer $M$ such that $\max_i|u_i^m|\leq 1/2$ for $m\geq M$.
By the Taylor series expansion with remainder, we have
$\log (1+u)=u-u^2/2+3\theta u^3$
for $|u|\leq 1/2$, where $|\theta |=|\theta (u)|\leq 1$.
Then we get the bound
$$
\Big |\log \Big (\prod_{i=1}^{l(m)}(1+u_i^m)\Big )
-[g(t)-g(0)-\frac 12s_2(g;\lambda (m))]\Big |
$$
$$
\leq\sum_{i=1}^{l(m)}\big |\log (1+u_i^m)-u_i^m+\frac 12 (u_i^m)^2
\big |\leq 3\sum_{i=1}^{l(m)}|u_i^m|^3\leq 3v_p(g)
\max_i|u_i^m|^{3-p}\tag2.17
$$
for all $m\geq M$.
This yields statement (I) because $g$ is continuous
and $g\in\Cal W_p([0,T])$. 

To prove statement (II), suppose that $b:=b_{\lambda}(g)$
is defined and continuous on $[0,T]$.
Given $t\in (0,T]$, for each $m\geq 1$, let $l(m)$ be as before,
$P:=\Cal E_{\lambda}(g)$ and let $v_i^m:=[P(t_i^m\wedge t)-P(t_{i-1}^m)]
/P(t_{i-1}^m)$ for $i=1,\dots,l(m)$.
Since $g$ and $b$ are continuous, there exists an 
integer $M$ such that $\max_i|v_i^m|\leq 1/2$ for $m\geq M$.
As in (2.17), since $|e^x-1|\leq|x|e^{|x|}$ for $x\in\RR$, we get
$$
\Big |\log P(t)-\sum_{i=1}^{l(m)}v_i^m+\frac 
12\sum_{i=1}^{l(m)}[v_i^m]^2\Big |\leq 3\sum_{i=1}^{l(m)}|v_i^m|^3
$$
$$
\leq Cv_p(g;[0,T])\max_i|g(t_i^m\wedge t)-g(t_{i-1}^m)|^{3-p}
+Cb(T)\max_i|b(t_i^m\wedge t)-b(t_{i-1}^m)|^2
$$
for some constant $C$ and all $m\geq M$.
Since the right side tends to zero with $m\to\infty$, the limit (2.4)
exists for $f=P$ because
$$
\lim_{m\to\infty}\sum_{i=1}^{l(m)}[v_i^m]^2=\lim_{m\to\infty}
\sum_{i=1}^{l(m)}\big [g(t_i^m\wedge t)-g(t_{i-1}^m)-\frac 12
[b(t_i^m\wedge t)-b(t_{i-1}^m)]\big ]^2=b(t).
$$
It then follows that
$$
\Cal L_{\lambda}\big (\Cal E_{\lambda}(g)\big )(t)=\lim_{m\to\infty}
\sum_{i=1}^{l(m)}v_i^m=\log P(t)+\frac 12\lim_{m\to\infty}
\sum_{i=1}^{l(m)}[v_i^m]^2=g(t)-g(0).
$$
The proof of Proposition 2.8 is complete.
\qed\enddemo

{\sl 2.4. Price and return}.
We define a price and its return in a duality under minimal
restrictions on stochastic processes.
To begin with we define a random moment $\tau_P$
which can be interpreted as the time of the crash  of a stock.
In the case when a stock price $P$ is the solution to the Dol\'eans--Dade 
equation (1.1), $\tau_P$ is the first moment $t$ when $P(t)\leq 0$.
Given a stochastic process $X=\{X(t)\colon\,t\in [0,T]\}$, let
$$
\tau_P=\tau_P(X):=
\cases t,&\text{if $\inf_{s\in [0,t)}X(s)>0$ and $X(t)\leq 0$ 
for $t\in (0,T]$,}\\t+,&\text{if $\inf_{s\in [0,t]}X(s)>0$ and 
$X^{-}(t+)\leq 0$ for $t\in (0,T)$,}\\
T+,&\text{if $\inf_{s\in [0,T]}X(s)>0$,}\endcases\tag2.18
$$
where $f^{-}(t+):=\lim\sup_{s\downarrow t}f(s)$.
Let $K$ be the set consisting of points $0$, $0+$, $t-$, $t$, $t+$ 
for $t\in (0,T]$ with the natural linear ordering: $s+<t-<t<t+$ if $s<t$, 
and endowed with the interval topology. 
Then $\tau_P$ is the random variable with values in $K$.
If $X$ is a price, then the event $\{\tau_P=t\}$ can be interpreted as 
the crash right before the time $t$, while the event $\{\tau_P=t+\}$ 
can be interpreted as the crash right after the time $t$.
Also, let $[0,t+):=[0,t]$.

Recalling notation $\Q([0,T])$, $\Cal L_{\lambda}$ and $\Cal E_{\lambda}$
from Definitions 2.2 and 2.5, we have:

\definition{Definition 2.9}
Let $R=\{R(t)\colon\,t\in [0,T]\}$ and $\spr=\{\spr (t)\colon\,t\in [0,T]\}$ 
be stochastic processes on a complete probability space 
$(\Omega,\Cal F,\Pr )$ such that $R(0)=0$ and $P(0)=1$ almost surely.
\roster
\item The pair $(P,R)$ will be called the {\it weak evolutionary
system} on $[0,\tau_P)$ if, for each $\lambda\in\Q([0,T])$, there exists 
$N=N(\lambda)\in\Cal F$ with $\Pr (N)=0$ such that,
for each $\omega\in\Omega\setminus N$,
the functions $\Cal E_{\lambda}(R(\cdot,\omega))$, 
$\Cal L_{\lambda}(P(\cdot,\omega))$ are defined on $[0,T]\cap 
[0,\tau_P(P(\omega)))$, and  satisfy the relations
$$
P(t,\omega)=\Cal E_{\lambda}(R(\cdot,\omega))(t)\quad\text{and}\quad
R(t,\omega)=\Cal L_{\lambda}(P(\cdot,\omega))(t)
$$
for each $t\in [0,T]\cap [0,\tau_P(P(\omega)))$.
\item If, in addition, the null set $N\in\Cal F$ in {\rm (1)} can be 
chosen the same for all $\lambda\in\Q([0,T])$, then the pair $(P,R)$ will 
be called the {\it evolutionary system} on $[0,\tau_P)$.
\endroster
If the pair $(P,R)$ is the weak evolutionary system then we call
$P$ the price and $R$ the return.
\enddefinition

We show that if $R$ is a Brownian motion and $P$ is a geometric
Brownian motion, then the pair $(P,R)$ is the weak evolutionary
system but not the evolutionary system (Proposition 2.11
and Remark 2.12).
However, according to the following statement, 
if almost all sample {\paths} of $R$ are in $\Cal W_2^{\ast}([0,T])
\cap\Cal D([0,T])$ then $(P,R)$ is the evolutionary system.
Recall that $f\in\Cal D([0,T])$ if, at each point of $(0,T)$,
$f$ is is either right--continuous or left--continuous.

Next we define a random moment $\tau_R$ for a return.
Given a stochastic process $Y=\{Y(t)\colon\,t\in [0,T]\}$
with almost all sample {\paths} in $\Cal D([0,T])$, let
$$
\tau_R=\tau_R(Y):=
\cases t,&\text{if $\inf_{s\in [0,t)}\Delta Y(s)>-1$ 
and $\Delta^{-}Y(t)\leq -1$ for $t\in (0,T]$,}\\
t+,&\text{if $\inf_{s\in [0,t)}\Delta Y(s)>-1$ and 
$\Delta^{+}Y(s)\leq -1$ for $t\in (0,T)$,}\\
T+,&\text{if $\inf_{s\in [0,T]}\Delta Y(s)>-1$.}\endcases
$$
Here $\Delta Y(s):=0$ everywhere except at jump points $s$ 
of $Y$ where either $\Delta Y(s):=\Delta^{+}Y(s)$ if it is non-zero,
or $\Delta Y(s):=\Delta^{-}Y(s)$ if it is non-zero.

\proclaim{Proposition 2.10}
Let $R=\{R(t)\colon\,t\in [0,T]\}$ be a stochastic process on
a complete probability space $(\Omega,\Cal F,\Pr)$ such that
$R(0)=0$ and  $R\in\Cal W_2^{\ast}([0,T])\cap\Cal D([0,T])$ with
probability $1$.
Then the indefinite product integral $\Cal P(R)$ is defined
with respect to almost every sample function of $R$,
$\tau_P(\Cal P(R))=\tau_R(R)$ almost surely and the pair
$(\Cal P(R),R)$ is the evolutionary system on $[0,\tau_R)$.
\endproclaim

\demo{Proof}
Let $N\in\Cal F$ be such that $\Pr (N)=0$ and $R(\cdot,\omega)
\in\Cal W_2^{\ast}([0,T])\cap\Cal D([0,T])$ for all 
$\omega\in\Omega\setminus N$.
Let $\lambda=\{\lambda (m)\colon\,m\geq 1\}\in\Q([0,T])$.
For each $t<\tau (\omega):=\tau_R(R(\omega))$, let $P(t,\omega):=
\Cal E_{\lambda}(R(\cdot,\omega))(t)$
if $\omega\in\Omega\setminus N$ and $P(t,\omega ):=0$ if $\omega\in N$.
By Proposition 2.6, $P(t,\omega )=\prodi_0^t(1+dR(\cdot,\omega))$ for 
$t<\tau (\omega)$ and $\omega\in\Omega\setminus N$.
For all $\omega\in\Omega$ and $t\geq\tau (\omega)$, let
$P(t,\omega):=\prodi_0^{\tau (\omega)}(1+dR(\cdot,\omega))$. 
Then $P=\{P(t)\colon\,t\in [0,T]\}$ is a stochastic process
by construction.
Fix $\omega\in\Omega\setminus N$ and let $\tau:=\tau (\omega)$, 
$P(t):=P(t,\omega)$, $R(t):=R(t,\omega)$.
Then $P(\tau -)$ exists and
$\inf\{P(t)\colon\,t\in [0,\tau )\}\geq\delta$ for some $\delta>0$
by Proposition 4.30 of Dudley and Norvai\v sa (1999).
By Lemmas 5.1 and 5.2 of Dudley and Norvai\v sa (1999), 
$P$ and $R$ have the same jump points on $[0,\tau)$ and thus 
$P\in\Cal D([0,\tau))$.
Moreover, $P\in\Cal W_2^{\ast}([0,\tau ))$ by Theorem 2.7.
By Proposition 2.3, $\widetilde R(t):=\Cal L_{\lambda}(P)(t)$ is defined
for each $t<\tau$ and has the same value with the indefinite
$LY$ integral $\Cal S(P)(t)$.
By Theorem 2.7 again, $\widetilde R(t)
=(LY)\smallint_0^tP^{-1}\,dP=R(t)$ for each $t<\tau$.
Since the null set $N$ does not depend on $\lambda$, 
the pair $(P,R)$ is evolutionary system on $[0,\tau_P)$.
It is clear that $\tau_P(P(\omega))=\tau_R(R(\omega))$ for all
$\omega\in\Omega\setminus N$.
The proof of Proposition 2.10 is complete.
\qed\enddemo

For example, a fractional Brownian motion $B_H$ with $H\in (1/2,1)$ 
and a symmetric $\alpha$-stable L\'evy motion $X_{\alpha}$ with 
$\alpha\in (0,2)$ are the returns of the evolutionary systems
$(P_H,B_H)$ and $(P_{\alpha},X_{\alpha})$, respectively, where
$$
P_H(t):=\exp\{B_H(t)\}\quad\text{and}\quad P_{\alpha}(t):=\exp\{
X_{\alpha}(t)\}
\prod_{(0,t]}(1+\Delta^{-}X_{\alpha})\exp\{-\Delta^{-}X_{\alpha}\}
$$
for $t\in [0,T]$.
The price $P_{\alpha}$ is positive until the first
moment $t$ when $\Delta^{-}X_{\alpha}(t)\leq -1$.

\proclaim{Proposition 2.11}
Let $B=\{B(t)\colon\,t\in [0,T]\}$ be a standard Brownian motion
and let $P_B:=\{\exp\{B(t)-t/2\}\colon\,t\in [0,T]\}$.
Then the pair $(P_B,B)$ is the weak evolutionary system on $[0,T]$.
\endproclaim

\remark{Remark {\rm 2.12}}
The proof of the above proposition rely on Th\'eor\`eme 5 of L\'evy 
(1940, p.\ 510):
for a standard Brownian motion $B$ and for a sequence $\lambda
=\{\lambda (m)\colon\,m\geq 1\}\in\Q([0,1])$, the limit
$\lim_{m\to\infty}s_2(B;\lambda (m))=1$ exists with probability 1.
However, the exceptional null set $N(\lambda)\in\Cal F$ of
this implication depends on $\lambda$ and $\cup\{N(\lambda)\colon\,
\lambda\in\Q([0,1])\}=\Omega$.
Moreover, for almost all $\omega\in\Omega$ there exist $\lambda\in
\Q([0,1])$ such that $\lim_{m\to\infty}s_2(B(\cdot,\omega);\lambda(m))
=\infty$, and hence, $\Cal E_{\lambda}(B(\cdot,\omega))=0$.
The proofs of these properties are given by Freedman (1983, p. 48)
because $\cup_m\lambda (m)$ is everywhere dense in $[0,1]$ if and
only if the mesh $|\lambda (m)|\to 0$ with $m\to\infty$.
\endremark

\demo{Proof}
The claim will follow from Proposition 2.8 once
we show that, given $\lambda=\{\lambda (m)\colon\,m\geq 1\}\in\Q([0,T])$,
the limit (2.14) with $g=B$ exists with probability 1 for each $t\in [0,T]$.
For each $t\in (0,T]$ and $m\geq 1$, let
$b_m(t,\omega):=\sum_{i=1}^{n(m)}\big [B(t_i^m\wedge t,\omega)-B(t_{i-1}^m
\wedge t,\omega)\big ]^2$.
Let $N=N(\lambda)\in\Cal F$ be a null set such that, for each $\omega\in
\Omega\setminus N$, $B(\cdot,\omega)$ is continuous function of bounded
$p$-variation for some $2<p<3$, and  the limit 
$\lim_mb_m(t,\omega)=t$ exists for each $t$ in the countable set 
$S:=\cup_m\lambda (m)$.
Let $t\in (0,T)\setminus S$.
Since $b_m(t,\omega)-t=[B(t,\omega)-B(t_{i-1}^m,\omega)]^2
+b_m(t_{i-1}^m,\omega)-t_{i-1}^m+[t_{i-1}^m-t]$ for $t_{i-1}^m
<t<t_i^m$, it follows that the limit as 
$m\to\infty$ of $b_m(t,\omega)$ is $t$ for each $t\in [0,T]$.
Therefore an appeal to Proposition 2.8 completes the proof.
\qed\enddemo

\vskip 1pc
\flushpar
{\bf 3 Pathwise trading strategies}
\vskip 1pc

\flushpar
A trading strategy is a collection of instructions for buying
and selling a stock, depending on its price fluctuations.
A mathematical notion of a trading strategy should be defined
so that one can calculate the portfolio value and portfolio gain for each
single trajectory of a stock price.
In the stochastic exponential model, the portfolio gain is the stochastic 
integral of a trading strategy with respect to a price.
Its value can be approximated by portfolio gains based on simple
trading strategies in probability.
In this section we take a pathwise approach to
trading strategies.

Consider a frictionless stock market with $\dim +1$
non-dividend-paying stocks and open for trade during a time 
period $[0,T]$.
A vector function $\spr=(\spr_0,\dots,\spr_{\dim})$ defined on $[0,T]$
will be called the {\it price} during the time period $[0,T]$
if, for each $k=0,\dots,\dim$, $\inf\{\spr_k(t)\colon\,t\in [0,T]\}>0$.
The value $\spr_k(t)$ refers to the price of the $k$th stock at time 
$t\in [0,T]$ for $k=0,\dots,\dim$.
For example, the price may be a vector of sample {\paths} of
the price stochastic processes defined in the preceding section.
A possible dependence of the price on a randomness and the duality
relation between the price and its return are
disregarded in this section.

As before, $\Q([0,T])$ denotes the set of all nested sequences 
$\lambda=\{\lambda (m)\colon\,m\geq 1\}$ of partitions 
of $[0,T]$ such that $\cup_m\lambda (m)$ is dense in $[0,T]$.
A sequence $\kappa=\{\kappa (m)\colon\,m\geq 1\}\in\Q([0,T])$
is a {\it refinement} of a sequence $\lambda=\{\lambda (m)\colon\,
m\geq 1\}\in\Q([0,T])$ if each $\kappa (m)$ is a refinement
of $\lambda (m)$.

\definition{Definition 3.1}
Let $\spr=(\spr_0,\dots,\spr_{\dim})$ be a price during a time 
period $[0,T]$.
Given $\lambda\in\Q([0,T])$, a vector function 
$\trs=(\trs_0,\dots,\trs_{\dim})$ defined on $[0,T]$
will be called the {\it $(\lambda,\spr)$--trading strategy}  
during the time period $[0,T]$ if, 
for each $k=0,\dots,\dim$ and $t\in [0,T]$, there exists the limit
$$
(LCS)\int_0^t\trs_k\,d_{\lambda}\spr_k:=\lim_{m\to\infty}\sum_{i=1}^{n(m)}
\trs_k(t_{i-1}^m\wedge t)[\spr_k(t_i^m\wedge t)-\spr_k(t_{i-1}^m
\wedge t)],\tag3.1
$$
where $\lambda=\{\lambda (m)\colon\,m\geq 1\}$ and $\lambda (m)
=\{t_i^m\colon\,i=0,\dots,n(m)\}$.
If there exists $\lambda_0\in\Q([0,T])$ such that 
for each refinement $\lambda\in\Q([0,T])$ of $\lambda_0$, $\phi$
is the $(\lambda,P)$--trading strategy on $[0,T]$ and (3.1) does
not depend on $\lambda$, then  
we call $\trs$ the {\it $\spr$--trading strategy}
during the time period $[0,T]$ and replace $d_{\lambda}$
with $d$ in the left side of (3.1).
\enddefinition

F\"ollmer (1981) proved that (3.1) exists whenever 
$\trs_k=f{\circ}\spr_k$ for some $f\in C^1$ and the quadratic variation 
is defined for $\spr_k$ along the sequence $\lambda\in\Q ([0,T])$.
The above notion of $(\lambda,\spr)$--trading strategy is similar
to the convergence of trading strategies introduced by
Bick and Willinger (1994, p. 356).
These authors derived the Black and Scholes 
formula without probabilistic arguments using
F\"ollmer's variant of It\^o's formula.
Remark 2.12 ensure that $\trs$ may be $(\lambda,\spr)$--trading
strategy without being a $\spr$-trading strategy.
It is clear that each $\trs_k$ is Left Cauchy--Stieltjes
integrable with respect to $\spr_k$ if $\trs$ is $\spr$--trading
strategy.
As the rest of this section show we could use the $LY$ integral instead
of the $LCS$ integral in the definition of $\spr$--trading strategies.
However in econometric analysis it seems easier to handle  with the
latter integral because to evaluate the $LY$ integral we would
need to know jumps of the price.
The following statement provides sufficient conditions for
a vector function to be the $\spr$--trading strategy.

\proclaim{Proposition 3.2} 
Let $\spr=(\spr_0,\dots,\spr_{\dim})$ be a price during
a time period $[0,T]$.
A vector function $\trs=(\trs_0,\dots,\trs_{\dim})$ on $[0,T]$ is
the $\spr$--trading strategy, each $\trs_k$ is $LY$ integrable with 
respect to $\spr_k$ on $[0,T]$ and 
$$
(LCS)\int_0^t\trs_k\,d\spr_k =(LY)\int_0^t\trs_k\,d\spr_k\tag3.2
$$
for each $t\in [0,T]$ and each $k=0,\dots,\dim$ in either of the 
following two cases{\rm :}
\roster
\item for each $k=0,\dots,\dim$,
$\spr_k\in\Cal W_2^{\ast}([0,T])$ and
$\trs_k=f_k{\circ}\spr_k$ for some $f_k\colon\,
\RR\mapsto\RR$ satisfying the local Lipschitz condition.
\item for each $k=0,\dots,\dim$, $\spr_k\in\Cal W_p([0,T])$ and
$\trs_k\in\Cal W_q([0,T])$ with $p,q >0$, $1/p+1/q>1$.
\endroster
\endproclaim

\demo{Proof}
In case (2) the conclusion follows from Theorem 2 and Corollary 3
of Norvai\v sa (1999).
To prove the conclusion in case (1),
for notation simplicity we suppress the index $k=0,\dots,\dim$
for $f_k$ and $\spr_k$.
By Theorem 2.1, $f{\circ}\spr$ is $LY$ integrable with respect to 
$\spr$ on $[0,T]$ and
$$
(LY)\int_0^t(f{\circ}\spr)\,d\spr=F{\circ}\spr (t)-F{\circ}P(0)-
\sum_{(0,t]}\phi^{-}-\sum_{[0,t)}\phi^{+},
$$
where $\phi^{-}:=\Delta^{-}(F{\circ}\spr)-(f{\circ}\spr)_{-}\Delta^{-}
\spr$, $\phi^{+}:=\Delta^{+}(F{\circ}\spr)-(f{\circ}\spr)\Delta^{+}\spr$
and $F(u):=\smallint_0^uf(x)\,dx$ for $u\geq 0$.
Fix $t\in (0,T]$ and let $\epsilon >0$.
By Lemma 2.4, there exists $\lambda=\{s_j\colon\,j=0,\dots,m\}\in
Q([0,t])$ such that $\sum_{j=1}^mv_2(P;(s_{j-1},s_j))<\epsilon$, and
for each refinement $\{t_i\colon\,i=0,\dots,n\}$ of $\lambda$,
$$
\Big |\sum_{(0,t]}\phi^{-}+\sum_{[0,t)}\phi^{+}
-\sum_{i=1}^{n}\big [\phi^{-}(t_i)+\phi^{+}(t_{i-1})\big ]\Big |
<\epsilon.
$$
Then choose $\{u_{j-1},v_j\colon\,j=1,\dots,m\}\subset [0,t]$
such that $s_{j-1}<u_{j-1}<v_j<s_j$ for $j=1,\dots,m$ and
$$
\sum_{j=1}^m\Big [\osc\,(P;(z_{j-1},u_{j-1}])+\osc\,(P;[v_j,z_j))\Big ]
<\epsilon.
$$
Then by the mean value theorem and using the Lipschitz condition
with the constant $K$, we have
$$
\Big |S_{LCS}(f{\circ}\spr,\spr;\kappa)-\big [F{\circ}\spr (t)
-F{\circ}\spr (0)-\sum_{(0,t]}\phi^{-}-\sum_{[0,t)}\phi^{+}\big ]\Big |
$$
$$
<\epsilon +\sum_{i=1}^n\Big |f{\circ}\spr(t_{i-1})[\spr (t_i)-
\spr (t_{i-1})] -[F{\circ}\spr (t_i)-F{\circ}\spr (t_{i-1})]
-[\phi^{-}(t_i)-\phi^{+}(t_{i-1})]\Big |
$$
$$
=\epsilon +\sum_{i=1}^n\Big |
[f(\theta_i)-f(P(t_i-))][\spr (t_i-)-\spr (t_{i-1}+)]
+[f{\circ}\spr (t_i-)-f{\circ}\spr (t_{i-1})]
$$
$$
\times [\spr (t_i)-\spr (t_{i-1}+)]\Big |
<\epsilon +K\epsilon +2\sup_t|P(t)|K\epsilon +K\epsilon,
$$
where $\theta_i\in [\spr (t_{i-1}+)\wedge\spr (t_i-),
\spr (t_{i-1}+)\vee\spr (t_i-)]$.
Since $\epsilon$ is arbitrary the proof of Proposition 3.2 is complete.
\qed\enddemo

Having defined the  $\spr$--trading strategies via the Left 
Cauchy--Stieltjes integral, the following definition of self--financing 
strategy corresponds naturally to its counterpart in the stochastic
exponent model.

\definition{Definition 3.3}
Let $\spr=(\spr_0,\dots,\spr_{\dim})$ be a price 
during a time period $[0,T]$ and let 
$\trs =(\trs_0,\dots,\trs_{\dim})$ be the $\spr$--trading strategy. 
\roster
\item The real-valued functions $\prf^{\trs}=\prf^{\trs,\spr}$ and 
$G^{\trs}=G^{\trs,\spr}$ defined on $[0,T]$ by
$$
\prf^{\trs}(t):=\sum_{k=0}^{\dim}\trs_k(t)\spr_k(t)\quad\text{and}
\quad G^{\trs}(t):=\sum_{k=0}^{\dim}(LCS)\int_0^t\trs_k\,d\spr_k
$$
are called the {\it portfolio value function} and the 
{\it portfolio gain function}, respectively.
\item
The $\spr$--trading strategy $\trs$ is called {\it self--financing} if 
$\prf^{\trs}(t)=\prf^{\trs}(0)+G^{\trs}(t)$ for each $t\in [0,T]$.
\endroster
\enddefinition

For the sake of illustration, next we partially extend Proposition 3.24 
of Harrison and Pliska (1981) to the present setting.
Let us denote the discounted price  
$(1,\spr_1/\spr_0,\dots,\spr_{\dim}/\spr_0)$ by $\overline\spr$.

\proclaim{Proposition 3.4}
Let $0<p<2$ and let $\spr =(\spr_0,\dots,\spr_{\dim})$ be a price
during a time period $[0,T]$ such that $\spr_k\in\Cal W_p([0,T])$
for $k=0,\dots,\dim$.
Suppose that $\trs=(\trs_0,\dots,\trs_{\dim})$ is a vector function
on $[0,T]$ such that $\trs_k\in\Cal W_p([0,T])$ for $k=0,\dots,\dim$.
Then $\trs$ is self--financing $\spr$--trading strategy if and only if 
it is self--financing $\overline\spr$--trading strategy.
\endproclaim

\demo{Proof}
Let $\beta :=1/\spr_0$.
We have $\beta\spr_k\in\Cal W_p([0,T])$ for each
$k=1,\dots,\dim$.
Thus, by case (2) of Proposition 3.2 with $q=p$, $\trs$ is the
$\spr$--trading strategy if and only if $\trs$ is the 
$\overline\spr$--trading strategy

To prove the ``only if'' part of the statement
suppose $\trs$ is self--financing $\spr$-trading strategy.
Let $V:=V^{\trs,\spr}$ and $\overline V:=V^{\trs,\overline\spr}
=\beta V$.
By Theorem 2.1 with $F(u)=u_1u_2$ for $u=(u_1,u_2)$, $h=
\trs_k$, and $\alpha =1$, for each $t\in [0,T]$, we have
$$
(LY)\int_0^t\trs_k\,d(\beta\spr_k)=(LY)\int_0^t\trs_k\beta\,
d\spr_k+(LY)\int_0^t\trs_k\spr_k\,d\beta
+\sum_{(0,t]}(\trs_k)_{-}\Delta^{-}\beta\Delta^{-}\spr_k\tag3.3
$$
for $k=0,\dots,\dim$, where the left side of (3.3) is equal to 0 when
$k=0$.
Since $\spr$--trading strategy $\trs$ is self--financing,
by the substitution rule for the Left Young integral (Theorem 9
of Norvai\v sa, 1999), for each $t\in [0,T]$, we have
$$
(LY)\int_0^t\beta\,dV=\sum_{k=0}^{\dim}(LY)\int_0^t\beta\,d\Big ((LY)
\int_0^{\cdot}\trs_k\,d\spr_k\Big )
=\sum_{k=0}^{\dim}(LY)\int_0^t\beta\trs_k\,d\spr_k.\tag3.4
$$
Using Theorem 2.1 again except that now $h=1$, and Proposition 7 
of Norvai\v sa (1999) about jumps of the indefinite $LY$ integral,
for each $t\in [0,T]$, we get
$$
\align
\overline V(t)-\overline V(0)
&=(LY)\int_0^t\beta\,dV+(LY)\int_0^tV\,d\beta
+\sum_{(0,t]}\Delta^{-}\beta\Delta^{-}V\\
\text{by (3.4)}\quad
&=\sum_{k=0}^{\dim}\Big\{(LY)\int_0^t\trs_k\beta\,d\spr_k+(LY)\int_0^t\trs_k
\spr_k\,d\beta
+\sum_{(0,t]}(\trs_k)_{-}\Delta^{-}\beta\Delta^{-}P_k\Big\}\\
\text{by (3.3)}\quad
&=\sum_{k=1}^{\dim}(LY)\int_0^t\trs_k\,d(\beta\spr_k)
=G^{\trs,\overline\spr}(t).
\endalign
$$
Thus $\trs$ is self--financing $\overline\spr$--trading strategy.
The proof of the converse implication is similar and therefore
is omitted.
\qed\enddemo

We finish with the main argument in favor of the pathwise approach
to trading strategies.
In their discussion of the notion of trading strategy, Harrison and Pliska 
(1981, Section 7) made several suggestions.
For example, it would be desirable to show that
a claim is attainable if and only if
it is the limit (in some appropriate sense) of claims generated
by simple self--financing strategies.
Duffie and Protter (1992) and Eberlein (1992) proved that
the portfolio gain processes are approximable
by their discrete counterparts under certain conditions.
Next we show a kind of approximation of a ``contingent claim"
by simple self--financing strategies in the present context.
A trading strategy $\trs =(\trs_0,\dots,\trs_{\dim})$ is {\it simple} if
each $\trs_k$ is a step function on $[0,T]$.
The idea of the following statement originated from
Harrison, Pitbladdo and Schaefer (1984, Proposition 9), where this 
claim is proved for price processes with continuous  
sample {\paths} of bounded variation.  

\proclaim{Theorem 3.5}
Let $\spr$ be a price during a time period $[0,T]$ 
and let $\trs$ be a vector function on $[0,T]$, both satisfying either
of the two conditions of Proposition {\rm 3.2}.
If $\inf\{\prf^{\trs}(t)\colon\,t\in [0,T]\}>0$ 
then there exists a sequence $\{\trs^N\colon\,N\geq 1\}$ of simple 
self--financing $\spr$--trading strategies such that $\prf^{\trs^N}(0)=
\prf^{\trs}(0)$ and $\lim_{N\to\infty}\prf^{\trs^N}(T)=V^{\trs}(T)$. 
\endproclaim

\demo{Proof} We start with the construction of the sequence of simple
self--financing $\spr$--trading strategies based on a given
nested sequence $\{\lambda^N\colon\,N\geq 1\}$ of partitions
$\lambda^N=\{0=t_0<\dots <t_n=T\}$.
Given an integer $N\geq 1$, we define $\trs^N=(\trs_0^N,\dots,
\trs_{\dim}^N)$ recursively with constant values on each interval of 
the partition $\lambda^N$.
For each $k=0,\dots,\dim$, let $\trs_k^N:=\trs_k(0)$ on $[0,t_1)$.
Suppose that all $\trs_k^N$ are defined on $[0,t_i)$ for some
$1\leq i\leq n$.
Let $\prf^N(t_i):=\sum_{k=0}^{\dim}\trs_k^N(t_{i-1})\spr_k(t_i)$.
Then for each $k=0,\dots,\dim$, let $\trs_k^N$ be equal to 
$\trs_k(t_i)\prf^N(t_i)/\prf^{\trs}(t_i)$ either on
$[t_i,t_{i+1})$ if $i<n$, or on $\{T\}$ if $i=n$.
It is clear that each $\trs^N$ is a simple $\spr$--trading strategy.
Moreover, the portfolio value function $\prf^{\trs^N}$ has values
$\prf^{\trs^N}(0)=\prf^{\trs}(0)$
and
$$
\prf^{\trs^N}(t_i)=\sum_{k=0}^{\dim}\trs_k^N(t_i)\spr_k(t_i)
=\frac {\prf^N(t_i)}{\prf^{\trs}(t_i)}\sum_{k=0}^{\dim}\trs_k(t_i)
\spr_k(t_i)=\prf^N(t_i)\tag3.5
$$
for each $i=1,\dots,n$.
Next we show that each $\spr$--trading strategy $\trs^N$ is 
self--financing.
Let $u=t_{i-1}$ and $v\in (t_{i-1},t_i]$ for some
$i=1,\dots,n$.
Since $\trs_k^N$ is constant on $[u,v)$ we get
$$
(LY)\int_u^v\trs_k^N\,d\spr_k=(RS)\int_u^v(\trs_k^N)_{-}^{(u)}\,
d(\spr_k)_{+}^{(v)}+\trs_k^N(u)\Delta^{+}\spr_k(u)
$$
$$
=\trs_k^N(u)[\spr_k(v)-\spr_k(u+)]+\trs_k^N(u)\Delta^{+}\spr_k(u)
=\trs_k^N(u)[\spr_k(v)-\spr_k(u)]
$$
for each $k=0,\dots,\dim$.
Given $t\in (0,T]$, let $l:=\max\{i\leq n\colon\,t_i\leq t\}$.
Then using the additivity of the $LY$ integral over adjacent intervals
(Theorem 4 of Norvai\v sa, 1999) and changing the 
order of summation over $k$ and $i$, we get
$$
\align
G^{\trs^N}(t)&=\sum_{k=0}^{\dim}\Big\{(LY)\int_{t_{l}}^t\trs_k^N\,d\spr_k
+\sum_{i=1}^{l}(LY)\int_{t_{i-1}}^{t_i}\trs_k^N\,d\spr_k\Big\}\\
&=\sum_{k=0}^{\dim}\trs_k^N(t_{l})[\spr_k(t)-\spr_k(t_{l})]
+\sum_{i=1}^{l}\big [V^N(t_i)-V^{\trs^N}(t_{i-1})\big ]\\\
\text{by (3.5)}\quad
&=\sum_{k=0}^{\dim}\big [\trs_k^N(t_{i(t)})\spr_k(t)-\trs_k^N(0)\spr_k(0)
\big ]=\prf^{\trs^N}(t)-\prf^{\trs^N}(0).
\endalign
$$
Thus the $\spr$--trading strategy $\trs^N$ is self--financing for each 
$N\geq 1$.
By (3.5) and by the additivity of the $LY$ integral over 
adjacent intervals again, we get  
$$
\Delta^N (t_{i-1},t_i):=\prf^{\trs} (t_i)-\prf^{\trs} (t_{i-1})
\frac {\prf^{\trs^N}(t_i)}{\prf^{\trs^N}(t_{i-1})}
=\prf^{\trs}(t_i)-\sum_{k=0}^{\dim}\trs_k(t_{i-1})\spr_k(t_i)
$$
$$
=\sum_{k=0}^{\dim}\Big\{(LY)\int_{t_{i-1}}^{t_i}\trs_k\,d\spr_k
-\trs_k(t_{i-1})[\spr_k(t_i)-\spr_k(t_{i-1})]\Big\}
=:\sum_{k=0}^{\dim}\Delta_k^N(t_{i-1},t_i)\tag3.6
$$
for each $i=1,\dots,n$.
Suppose that one can choose a sequence $\{\lambda^N\colon\,N\geq 1\}$
such that
$$
\max_{1\leq i\leq n}|\Delta^N (t_{i-1},t_i)|\leq\sum_{i=1}^{n}
|\Delta^N (t_{i-1},t_i)|\leq \epsilon_N\tag3.7
$$
for some $\epsilon_N\downarrow 0$.
Then it follows that
$$
\Big |\frac {\prf^{\trs^N}(t_i)}{\prf^{\trs}(t_i)}\frac {\prf^{\trs}
(t_{i-1})}{\prf^{\trs^N}(t_{i-1})}-1\Big |
=\Big |\frac {\Delta^N (t_{i-1},t_i)}{\prf^{\trs}(t_i)}\Big |
\leq\epsilon_N/\delta
$$
for each $t_{i-1}, t_i\in\lambda^N$.
We conclude then recursively that $V^N(t_i)\neq 0$
for $i=1,\dots,n$ whenever $2\epsilon_N\leq\delta$.
By the mean value theorem, $|\log (1+u)|\leq 2|u|$ for each
$|u|\leq 1/2$.
Thus, for all $N$ such that $2\epsilon_N\leq\delta$, using the 
telescoping sum representation, we get
$$
\align
\Big |\log\frac {\prf^{\trs^N}(T)}{\prf^{\trs}(T)}\Big |
&\leq\sum_{i=1}^{n}\Big |\log\frac {\prf^{\trs^N}(t_i)}{\prf^{\trs}
(t_i)}-\log\frac {\prf^{\trs^N}(t_{i-1})}{\prf^{\trs}(t_{i-1})}\Big|
\leq 2\sum_{i=1}^{n}\Big |\frac {\prf^{\trs^N}(t_i)}{\prf^{\trs}(t_i)}
\frac {\prf^{\trs}(t_{i-1})}{\prf^{\trs^N}(t_{i-1})}-1\Big |\\
&\leq\frac {2}{\delta}\sum_{i=1}^{n}
|\Delta^N (t_{i-1},t_i)|\leq 2\epsilon_N/\delta,
\endalign
$$
where the last inequality follows from the second inequality in (3.7).
Since $\prf^{\trs^N}(0)=\prf^{\trs}(0)$ by construction,
this yields the conclusion of the theorem.

It remains to find $\{\lambda^N\colon\,N\geq 1\}$ such that (3.7)
holds for some $\epsilon_N\downarrow 0$.
Suppose that condition (1) of Proposition 3.2 holds.
By Theorem 2.1 and by the mean value theorem, each term $\Delta_k^N
(u,v)$ in (3.6) with $u=t_{i-1}$, $v=t_i$, $i=1,\dots,n$ and 
$k=0,\dots,\dim$ is equal to
$$
[f_k(\theta_k)-f_k(\spr_k(v-))][\spr_k(v-)-\spr_k(u+)]
+[f_k(\spr_k(v-))-f_k(\spr_k(u))][\spr_k(v)-\spr (u+)]
$$
$$
-\sum_{(u,v)}\big [\Delta^{-}(F_k{\circ}\spr_k)-(f_k{\circ}\spr_k)_{-}
\Delta^{-}\spr_k\big ]-\sum_{(u,v)}\big [\Delta^{+}(F_k{\circ}\spr_k)
-(f_k{\circ}\spr_k)\Delta^{+}\spr_k\big ],
$$
where $\theta_k\in (\spr_k(u+)\wedge\spr_k(v-),
\spr_k(u+)\vee\spr_k(v-)]$ and $F_k(u)=\smallint_0^uf_k(x)\,dx$
for $u\geq 0$.
Then, given $\epsilon_N\downarrow 0$, one can find $\lambda^N$ such that
(3.7) holds just as in the proof of Proposition 2.3.
Finally, suppose that condition (2) of Proposition 3.2 holds.
Choose $p'>p$ and $q'>q$ so that $1/p'+1/q'>1$.
By the Love--Young inequality (p.\ 256 in Young, 1936),
each term $\Delta_k^N(u,v)$ in (3.6)
with $u=t_{i-1}$, $v=t_i$, $i=1,\dots,n$ and $k=0,\dots,\dim$ 
can be bounded as follows:
$$
\align |\Delta_k^N(u,v)|
&=\Big |(RS)\int_u^v[(\trs_k)_{-}^{(u)}-(\trs_k)_{-}^{(u)}(u)]
\,d(\spr_k)_{+}^{(v)}+\sum_{(u,v)}\Delta^{-}\trs_k\Delta^{+}\spr_k\Big 
|\\&\leq KV_{p'}(\spr_k;[u+,v])V_{q'}(\trs_k;[u,v-])
+\sum_{(u,v)}|\Delta^{-}\trs_k\Delta^{+}\spr_k|
\endalign
$$
for some finite constant $K$ depending on $p'$ and $q'$ only.
Again, given $\epsilon_N\downarrow 0$, one can find $\lambda^N$ such that
(3.7) holds as in the proof of Proposition 2.3.
The proof of Theorem 3.5 is complete.
\qed\enddemo

\vskip 1pc
\flushpar
{\bf 4 Implications and conclusions}
\vskip 1pc

\flushpar
The results of the present paper provide an alternative construction
of a stock price model, and show that many concrete financial models 
can be treated using classical calculus.
By its definition, the evolutionary system is the continuous--time
model obtained as the limit of the discrete--time model (2.2)
along a sequence of partitions of a time interval into shrinking
subintervals.
The evolutionary system separates analytical and  probabilistic
aspects of analysis casting new light on important
problems of stock price modelling.
The arbitrage construction in Subsection 1.5 illustrates
implications of this separation.

New definition of the return $R$ in the evolutionary
system $(P,R)$ makes easier to use it in a statistical
analysis as compared with the definition (1.2).
Statistical analysis of analytical properties of functions
developed in relation to natural sciences could be applied
in econometric analysis of the evolutionary system.
For example, an interesting task is to distinguish 
the hypotheses that the $p$-variation index
$\pindex (R)<2$ against the hypotheses that $\pindex (R)\geq 2$.
This is important because the value $\pindex (R)=2$ separates
a fundamentally different behaviour of $R$.
Also, testing hypotheses $\pindex (R)<2$ and $R$ is continuous 
could be used to test market efficiency related to arbitrage.
Naturally that there are no ready to use statistical tests
for estimating the $p$-variation index.
In this case one needs to extract from data an information about a local
behavior of a sample function rather than an information about
tail distribution, or correlation estimates.
The first step in this direction has been taken up 
by Norvai\v sa and Salopek (1999). 
These authors suggest a statistic based on old results of G.\ Baxter 
and E.G.\ Gladyshev concerning quadratic variation for Gaussian processes.
Also, they compare the results of data analysis using
the new definition of the return and the log return.

Let dim$_{HB}(G)$ be  the Hausdorff--Besicovitch dimension of a set
$G$.
Then for a large class of stochastic processes,
the relation dim$_{HB}($graph $X)=2-1/(1\vee \pindex (X))$
holds for almost all sample functions of $X$.
This fact can be used to construct new statistics for
estimating the $p$-variation index $\pindex (X)$ because statistical 
analysis of fractal 
dimensions is already highly developed in various natural sciences.
The real analysis approach to modelling of stock price changes provides
a new meaning to stylized facts discovered in Econophysics
(see e.g.\ Bouchaud and Potters, 1999), and opens a way for exploring 
new tools for investigating financial markets.

\vskip 1pc
\flushpar
{\bf References}
\vskip 1pc

\noindent
\def\nref {\global\advance\q by1 \item{\bf\the\q.}}
\newcount\q \q=0

%\smallskip
%\hang{Alexander, S. S. (1961).
%Price Movements in Speculative Markets: Trends or Random Walks.
%{\sl Industrial Managemnet Review}, {\bf 2}(2), 7-26.} 

\smallskip
\nref Bick, A., Willinger, W.:
Dynamic spanning without probabilities.
Stoch. Proc. Appl. {\bf 50}, 349-374 (1994)

\smallskip
\nref Bouchaud, J.-P., Potters, M.:
Theory of financial risk: from data analysis to risk management.
Science \& Finance, 1999 (to appear)

\smallskip
\nref B\"uhlmann, H., Delbaen, F., Embrechts, P., Shiryaev, A.N.:
No--arbitrage, change of measure and conditional Esscher transforms.
CWI Quarterly {\bf 9}, 291-317 (1996)

\smallskip
\nref Campbell, J.Y., Lo, A.W., MacKinlay, A.C.:
The econometrics of financial markets.
Princeton, New Jersey: Princeton University Press 1997

\smallskip
\nref Clarkson, R. S.:
Financial economics - an investment actuary's viewpoint.
British Actuarial J. {\bf 2}, IV, 809-973 (1996)

\smallskip
\nref Clarkson, R. S.:
An actuarial theory of option pricing.
British Actuarial J. {\bf 3}, II, 321-410 (1997)

\smallskip
\nref Dol\'eans--Dade, C.:
Quelques applications de la formule de changement de variables pour
les semimartingales.
Z. Wahrsch. verw. Geb. {\bf 16}, 181-194 (1970)

\smallskip
\nref Dudley, R.M., Norvai\v sa, R.:
Product integrals, Young integrals and $p$-variation.
Lect. Notes in Math. {\bf 1703} Berlin: Springer 1999, pp 73-214

\smallskip
\nref Duffie, D., Protter, P.:
From discrete- to continuous-time finance: 
Weak convergence of the financial gain process.
Math.\ Finance {\bf 2}, 1-16 (1992)

\smallskip
\nref Eberlein, E.:
On modelling questions in security valuation.
Math.\ Finance {\bf 2}, 17-32 (1992)

\smallskip
\nref Elton E.J., Gruber M.J., Kleindorfer, P.R.:
A closer look at the implications of the stable paretian hypothesis.
Review of Economics and Statistics, 231-235 (1975)

\smallskip
\nref F\"ollmer, H.:
Calcul d'It\^o sans probabilites.
In: Az\'ema, J., Yor, M. (eds.); 
S\'eminaire de Probabilit\'es XV; Lect. Notes
in Math. {\bf 850}. Berlin: Springer 1981, pp 143-150

\smallskip
\nref Focardi, S., Jonas, C.:
Modeling the market: New theories and techniques.
New Hope, Pennsylvania: F.J. Fabozzi 1997

\smallskip
\nref Freedman, D.: Brownian motion and diffusion.
New-York: Springer 1983

\smallskip
\nref Harrison, J.M., Pliska, S.R.:
Martingales and stochastic integrals in the theory
of continuous trading.
Stoch. Proc. Appl. {\bf 11}, 215-260 (1981)

\smallskip
\nref Harrison, J.M., Pitbladdo, R., Schaefer, S.M.:
Continuous price processes in frictionless markets have
infinite variation.
Journal of Business {\bf 57}, 353-365 (1984)

\smallskip
\nref Kac, M., Rota, G.-C., Schwartz, J.T.:
Discrete thoughts. Essays on mathematics, science,
and philosophy.
Revised second edition.
Boston: Birkh\"auser 1992

\smallskip
\nref L\'evy, P.:
Le mouvement brownien plan.
Amer.\ J.\ Math. {\bf 62}, 487-550 (1940)

\smallskip
\nref Norvai\v sa, R.:
$p$-variation and integration of sample functions of stochastic 
processes.
In: Grigelionis, B.\  et al.\ (eds.); Prob. Theory and Math. Stat.,
1999 (to appear)

\smallskip
\nref Norvai\v sa, R., Salopek, D.M.:
Estimating the Orey index of a Gaussian stochastic process
with stationary increments: An application to financial data set.
In: Ivanoff, G.\ et al.\  (eds.); Proc.\ Int.\ Conf.\ on
Stochastic Models, 1999 (to appear)

\smallskip
\nref Pliska, S.R.:
Introduction to Mathematical Finance. Discrete time models.
USA: Blackwell Publishers 1997

\smallskip
\nref Salopek, D.M.:
Tolerance to arbitrage.
Stoch. Proc. Appl. 76, 217-230 (1998)

\smallskip
\nref Shiryaev, A.N.:
Essentials of stochastic finance.
Part 2. Theory.
Transl. from Russian by N. Kruzhilin.
World Scientific 1998 (to appear)

\smallskip
\nref Wong, E., Zakai, M.:
On the convergence of ordinary integrals to stochastic integrals.
Ann. Math. Statist. {\bf 36}, 1560-1564 (1965)

\smallskip
\nref Young, L.C.:
An inequality of the H\"older type, connected with Stieltjes
integration.
Acta Math. (Sweden) {\bf 67}, 251-282 (1936)

\enddocument